\documentclass[a4paper,11pt,leqno,twoside]{article}

\usepackage{amsmath}
\usepackage{amssymb}
\usepackage[margin=3.5cm]{geometry}

\usepackage{mathrsfs}
\renewcommand{\mathcal}[1]{\mathscr{#1}}

\newtheorem{lem}{Lemma}[section] \newtheorem{prp}[lem]{Proposition}
\newtheorem{df}[lem]{Definition} \newtheorem{thm}[lem]{Theorem} 
\newtheorem{crl}[lem]{Corollary}

\newlength{\eqdemoffset}
\newlength{\enumrqueoffset}
\newenvironment{dem}[1][] {\begin{trivlist}\item[]{\sc Proof#1.}~}
  {\ifhmode{\unskip~\nobreak}\else%
    {\nopagebreak\vspace{\eqdemoffset}\leavevmode\hfill}\fi%
    \hfill$\blacksquare$\end{trivlist}}
\newenvironment{rque}[1]{\refstepcounter{lem}%
  \begin{trivlist}\item[]{\sc #1 \arabic{section}.\arabic{lem}}~}%
  {\ifhmode{\unskip~\nobreak}\else%
    {\nopagebreak\vspace{\enumrqueoffset}\leavevmode\hfill}\fi%
    \hfill$\square$\end{trivlist}}

\renewcommand{\ss}{\scriptscriptstyle}
\newcommand{\ds}{\displaystyle}
\newcommand{\ts}{\textstyle}
\newcommand{\un}{1\!\!1}
\newcommand{\Ker}{\mathop{\mathrm{Ker}}}
\newcommand{\Img}{\mathop{\mathrm{Im}}}
\newcommand{\Irr}{\mathop{\mathrm{Irr}}}
\newcommand{\Dom}{\mathop{\mathrm{Dom}}}
\newcommand{\Der}{\mathop{\mathrm{Der}}}
\newcommand{\Inn}{\mathop{\mathrm{Inn}}}
\newcommand{\Supp}{\mathop{\mathrm{Supp}}}
\newcommand{\Tr}{\mathord{\mathrm{Tr}}}
\renewcommand{\iff}{\textbf{iff} } 
\newcommand{\jok}{\bigstar}
\newcommand{\red}{{\mathord{\mathrm{red}}}}
\newcommand{\Cst}{$C^*$-\relax}
\newcommand{\tens}{\mathord{\otimes}}
\newcommand{\id}{\mathord{\mathrm{id}}}
\newcommand{\rond}{\makebox[1em][c]{$\circ$}}

\renewcommand{\O}{E_1}
\newcommand{\B}{E_2}
\newcommand{\E}{E}
\newcommand{\uB}{{\underline E}_2}

\newcommand{\CC}{\mathbb C} 
\newcommand{\RR}{\mathbb R} 
\newcommand{\ZZ}{\mathbb Z}
\newcommand{\NN}{\mathbb N} 
\newcommand{\Aa}{\mathcal{A}}
\newcommand{\Cc}{\mathcal{C}}
\newcommand{\Dd}{\mathcal{D}}
\newcommand{\Hh}{\mathcal{H}}
\newcommand{\Ss}{\mathcal{S}}
\newcommand{\Kk}{\mathcal{K}}
\newcommand{\Vk}{\mathfrak{V}}
\newcommand{\Ek}{\mathfrak{E}}

\title{Paths in quantum Cayley trees \\ and $L^2$-cohomology}
\author{Roland Vergnioux} 
\date{}
\begin{document}

\maketitle

\begin{abstract}
  We study existence, uniqueness and triviality of path cocycles in the quantum
  Cayley graph of universal discrete quantum groups. In the orthogonal case we
  find that the unique path cocycle is trivial, in contrast with the case of
  free groups where it is proper. In the unitary case it is neither bounded nor
  proper. From this geometrical result we deduce the vanishing of the first
  $L^2$-Betti number of $A_o(I_n)$.

  MSC 2000: 46L65 (46L54, 20G42, 16E40, 20F65)
\end{abstract}

\smallskip

\section{Introduction}

Consider the unital $*$-algebras $\Ss = \Aa_o(I_n)$,
$\Aa_u(I_n)$ introduced by Wang \cite{Wang_Free_Prod} and defined by
$n^2$ generators $u_{ij}$ and relations as follows:
\begin{align*}
  \Aa_u(I_n) &= \langle u_{ij} ~|~ (u_{ij})_{ij}
  \text{~~and~~} (u_{ij}^*)_{ij} \text{~~are unitary}\rangle, \\
  \Aa_o(I_n) &= \Aa_u(I_n) / \langle u_{ij} = u_{ij}^*
  \text{~~for all~~} i,~ j\rangle.  
\end{align*}

These algebras naturally become Hopf algebras with coproduct given by
$\delta(u_{ij}) = \sum u_{ik}\tens u_{kj}$. It results from the general theory
of Woronowicz \cite{Woronowicz_Matrix} that they admit unique Haar integrals $h
: \Ss \to \CC$. We denote by $H$ the corresponding Hilbert space completions,
and by $S_\red$ the operator norm completions of the images of the natural
``regular'' representations $\lambda : \Ss \to B(H)$.

Since the introduction of $\Aa_o(I_n)$ and $\Aa_u(I_n)$ and the
seminal work of Banica \cite{Banica_Un,Banica_On} on their theory of
corepresentations, strong relations and similarities with the free group
algebras $\CC F_n$ have been noted:
\begin{itemize}
\item canceling the generators $u_{ij}$ with $i\neq j$ yields natural quotient
  maps $\Aa_u(I_n) \to \CC F_n$, $\Aa_o(I_n) \to \CC
  (\ZZ/2\ZZ)^{* n}$,
\item any Hopf algebra $\Ss$ associated with a unimodular ``discrete quantum
  group'' (resp.: and with self-adjoint generating corepresentation) is a Hopf
  quotient of some $\Aa_u(I_n)$ (resp. $\Aa_o(I_n)$),
\item the irreducible corepresentations of $\Aa_u(I_n)$ are naturally
  indexed by the elements of the free monoid on two generators $u$, $\bar u$.
\end{itemize}
Due to these remarks, the discrete quantum groups associated with
$\Aa_o(I_n)$ and $\Aa_u(I_n)$ have sometimes been called ``universal'', or even
``free''.

The analogy with free groups led in subsequent works to non-trivial \Cst
algebraic results showing that $A_o(I_n)_\red$ and $A_u(I_n)_\red$, $n\geq 2$,
indeed share important global analytical properties with free groups:
\begin{itemize}
\item the \Cst algebra $A_u(I_n)_\red$ is simple, the \Cst algebras
  $A_u(I_n)_\red$ and $A_o(I_n)_\red$ are non-nuclear (if $n\geq 3$)
  \cite{Banica_Un},
\item the discrete quantum groups associated with $\Aa_u(I_n)$,
  $\Aa_o(I_n)$ satisfy Akemann-Ostrand Property
  \cite{Vergnioux_Cayley},
\item the discrete quantum groups associated with $\Aa_u(I_n)$,
  $\Aa_o(I_n)$ have the Property of Rapid Decay \cite{Vergnioux_RD},
\item the \Cst algebra $A_o(I_n)_\red$ is simple (if $n\geq 3$) and exact, the
  von Neumann algebra $A_o(I_n)_\red''$ is a full and prime $II_1$ factor
  \cite{VaesVer_Boundary},
\item the von Neumann algebra $A_u(I_n)_\red''$ is a full and prime $II_1$
  factor \cite{VaesVennet_Boundary}.
\end{itemize}
The starting point of most results in \cite{Vergnioux_Cayley, VaesVer_Boundary,
  VaesVennet_Boundary} is the construction of quantum analogues of geometrical
objects associated with free groups, namely the quantum Cayley graphs and the
quantum boundaries of the discrete quantum groups associated with $\Aa_o(I_n)$,
$\Aa_u(I_n)$.

\bigskip

In the present paper, we pursue the study of the quantum Cayley graphs of
$\Aa_o(I_n)$ and $\Aa_u(I_n)$, concentrating in particular on paths to the
origin, which yield in good cases path cocycles, and we derive applications to
the $L^2$-cohomology of the algebras $\Aa_o(I_n)$, $\Aa_u(I_n)$. However, for
the first time the results we obtain for $\Aa_o(I_n)$ strongly contrast with the
known results for free groups. More precisely:
\begin{itemize}
\item on the geometrical side, we find out in Theorem~\ref{thm_path_Ao} that
  the path cocycle of $\Aa_o(I_n)$ is trivial, hence bounded, whereas the
  path cocyle on $F_n$ is unbounded, and even proper,
\item on the cohomological and analytical side, we prove in
  Theorem~\ref{thm_vanishing} that the first $L^2$-cohomology group
  $H^1(\Aa_o(I_n), M)$ vanishes, and in particular we have for the first
  $L^2$-Betti number $\beta_1^{(2)}(\Aa_o(I_n)) = 0$ for all $n$,
  whereas $\beta_1^{(2)}(F_n) = n-1$ and $\beta_1^{(2)}((\ZZ/2\ZZ)^{*n}) = \frac
  n2-1$.
\end{itemize}
In the case of $\Aa_u(I_n)$, we obtain a path cocycle which is not bounded, but
also not proper. We believe that these results open a new perspective in the
study of universal discrete quantum groups. In particular, they put a new
emphasis on the question of a-T-menability for $\Aa_o(I_n)$, which remains
open.\footnote{Note (January 2012). Since the writing of this article,
  a-T-menability of $\Aa_o(I_n)$ and $\Aa_u(I_n)$ has been established by
  Brannan~\cite{Brannan_aT}.}

Recall that in the classical case a-T-menability, also called Haagerup's
Property, has important applications, in particular in connection with
$K$-theory, and amounts to the existence of a metrically proper $\pi$-cocycle
for some unitary representation $\pi$ \cite{Book_aTmenability}. The results in
the present paper show that, in contrast with the case of free groups, such a
cocycle does not exist on $\Aa_o(I_n)$ if $\pi$ is contained in (a multiple of)
the regular representation.

Finally, let us mention a recent work of C.~Voigt \cite{Voigt_BC} about the
Baum-Connes conjecture and the K-theory of $A_o(I_n)$, which also presents the
feature that ``invariants do not depend on $n$'': namely it is shown that
$K_0(A_o(I_n)) = K_1(A_o(I_n)) = \ZZ$. It is also proved there that $A_o(I_n)$
is $K$-amenable, a new property shared with free groups.

\bigskip

Let us describe the strategy and the organization of the paper. After recalling
some notation about discrete quantum groups and classical Cayley graphs, we
introduce in Section~\ref{sec_cocycles} the notion of path cocycle for discrete
quantum groups and prove some elementary facts about bounded cocycles and path
cocycles. Then we continue the study of quantum Cayley trees started in
\cite{Vergnioux_Cayley}. We introduce in Section~\ref{sec_endpoints} a
``rotation operator'' and ``quasi-classical'' subspaces, and we give a new
description of the space of geometrical edges.

In Section~\ref{sec_cayley} we prove the existence and uniqueness of path
cocycles in quantum Cayley trees and we study their triviality. We prove in
particular the triviality of the path cocycle on $\Aa_o(I_n)$ for $n\geq 3$, and
we observe that the path cocycle on $\Aa_u(I_n)$ is unbounded. Finally a
``universality trick'' allows us to prove in Section~\ref{sec_cohomology} that
in fact the whole $L^2$-cohomology group $H^1(\Ss, M)$ vanishes in the
orthogonal case.

The reader interested in the application to $L^2$-cohomology will probably like
to skip the technical Section~\ref{sec_endpoints} at first, although it cannot
be completely avoided from a logical point of view. Let us also note that the
present article, except Section~\ref{sec_endpoints}, relies only on Sections~3
and~4 of \cite{Vergnioux_Cayley}. The results of this paper about $\Aa_o(I_n)$
have been announced at the Workshop on operator algebraic aspects of quantum
groups held in Leuven in November 2008. It is a pleasure to thank the organizer,
Stefaan Vaes.

\subsection{Notation}
\label{sec_notation}

Let us start from a regular multiplicative unitary $V \in B(H\tens H)$ with a
unique fixed line --- i.e. $V$ is of compact type without multiplicity
\cite{BaajSkand_Unitaires}. Our ``classical case'' will be the case of the
multiplicative unitary associated with a discrete group $\Gamma$ : then $H =
\ell^2\Gamma$ and $V = \sum_{g\in\Gamma} \un_g\tens\lambda_g$, where $\un_g$ is
the characteristic function of $\{g\}$ acting by pointwise multiplication on
$H$, and $\lambda_g$ is the left translation by $g$ in $H$. Note that all tensor
products between Hilbert spaces (resp. \Cst algebras, Hilbert \Cst modules) are
completed with respect to the hilbertian tensor norm (resp. the minimal \Cst
tensor norm).

A number of more familiar objects can be recovered from $V$. First of all we
have the Woronowicz \Cst algebra $(S_\red,\delta)$ \cite{Woronowicz_Matrix} and
the completion $(\hat S,\hat\delta)$ of the dual multiplier Hopf
algebra~\cite{VanDaele_Multiplier}:
\begin{align*}
  &S_\red = \overline{\mathrm{Span}} \{(\omega\tens\id)(V) ~|~ \omega\in
  B(H)_*\},
  ~~~ \delta(x) = V(x\tens 1)V^*; \\
  &\hat S = \overline{\mathrm{Span}} \{(\id\tens\omega)(V) ~|~ \omega\in
  B(H)_*\}, ~~~ \hat\delta(a) = V^*(1\tens a)V.
\end{align*}
In the classical case we recover $S_\red = C^*_\red(\Gamma)$ and $\hat S =
c_0(\Gamma)$ with their usual coproducts, hence the group $\Gamma$ itself. In
general we have $V \in M(\hat S\tens S_\red)$.

The theory of corepresentations of $(S_\red,\delta)$ or, equivalently, the \Cst
algebra structure of $\hat S$, is well understood via Peter-Weyl and
Tannaka-Krein type results \cite{Woronowicz_Matrix}. In particular, $\hat S$ is
a $c_0$-sum of matrix algebras and we denote by $\Cc$ the category of
finite-dimensional representations of $\hat S$. It is naturally endowed with a
monoidal structure for which every object $\alpha$ has a dual $\bar\alpha$. We
denote by $1_\Cc \in \Irr\Cc$ the trivial representation.

In the classical case, the irreducible representations of $\hat S$ have
dimension $1$ and correspond naturally to the elements of $\Gamma$ so that the
monoidal structure (resp. the duality operation, the trivial representation)
identifies to the group structure (resp. the inverse of the group, the unit
$e\in\Gamma$).

\bigskip

For $\alpha\in\Irr\Cc$ we denote by $p_\alpha \in \hat S$ the central support of
$\alpha$, so that $p_\alpha\hat S \simeq L(H_\alpha)$ and $\hat S =
c_0\text{-}\bigoplus_\alpha p_\alpha\hat S$. More generally to any finite subset
of $\Irr\Cc$ we associate a central projection of $\hat S$ by summing the
corresponding $p_\alpha$'s. In the classical case this projection is the
characteristic function of the considered subset of $\Irr\Cc \simeq \Gamma$. As
an exception we denote by $p_0$ the minimal central projection associated with
the trivial representation.

The algebraic direct sum of the matrix algebras $p_\alpha\hat S$ is a dense
subalgebra denoted $\hat \Ss$. Similarly the coefficients of
f.-d. corepresentations span a dense subalgebra $\Ss \subset S_\red$, and for
every $\alpha\in\Irr\Cc$ the product $(p_\alpha\tens\id)V$ lies in the algebraic
tensor product $\hat\Ss\tens\Ss$. In the classical case we have $\hat\Ss =
\CC^{(\Gamma)}$ and $\Ss = \CC\Gamma$.

We denote by $h$ the Haar state of $S_\red$. The Hilbert space $H$ identifies
with the GNS space of $h$ via a map $\Lambda_h : S_\red \to H$ given by the
choice of a fixed unit vector of $V$, $\xi_0 = \Lambda_h(1)$: we have then $h =
(\xi_0|\,\cdot\,\xi_0)$. We denote by $\Hh$ the dense subspace $\Lambda_h (\Ss)
\subset H$. In the classical case we can take $\xi_0 = \un_{e} \in
\ell^2\Gamma$.

From the rich structure of $(S_\red,\delta)$ we will also need the counit
$\epsilon : \Ss \to \CC$. In fact we will use $\epsilon$ at the hilbertian level,
hence we put $\epsilon(\Lambda_h(x)) = \epsilon(x)$ for $x \in \Ss$, thus
defining an unbounded linear form $\epsilon : \Hh \to \CC$. In the classical
case we have $\epsilon(\lambda_g) = 1$ and $\epsilon(\un_g) = 1$ for all
$g\in\Gamma$.

\bigskip

In this paper we are primarily interested in the case of the so-called universal
discrete quantum groups and their free products, see
\cite{Wang_Free_Prod,WangVanDaele}. More precisely, if $Q \in GL_n(\CC)$ is a
fixed invertible matrix, the unitary universal discrete quantum group $\Aa_u(Q)$
and the orthogonal universal discrete quantum group $\Aa_o(Q)$ are defined by
twisting the relations for $\Aa_o(I_n)$, $\Aa_u(I_n)$ given in the Introduction:
\begin{align*}
  \Aa_u(Q) &= \langle u_{ij} ~|~ (u_{ij})_{ij}
  \text{~~and~~} Q(u_{ij}^*)_{ij}Q^{-1} \text{~~are unitary}\rangle, \\
  \Aa_o(Q) &= \Aa_u(I_n) / \langle (u_{ij})_{ij} = Q(u_{ij}^*)_{ij}Q^{-1} \rangle.  
\end{align*}
In the case of $\Aa_o(Q)$ one requires moreover that $Q\bar Q$ is a scalar
multiple of $I_n$, and we will focus on the non-amenable cases, which correspond
to $n\geq 3$. These discrete quantum groups, and in particular their
corepresentation theory, have been studied in \cite{Banica_Un}.

\bigskip

We will recall definitions and facts about quantum Cayley graphs in the next
sections as they are needed : see {\sc Reminders}~\ref{rem_cayley},
\ref{rem_ascending}, \ref{rem_alg}, \ref{rem_Kg} and~\ref{rem_ql}. To give the
right perspective for these definitions, let us recall here some facts about the
classical case \cite{Serre_Arbres}.

A (classical) graph is given by a set of vertices $\Vk$, a set of edges $\Ek$,
an involutive reversing map $\theta : \Ek \to \Ek$ without fixed points, and an
endpoints map $e : \Ek \to \Vk^2$ such that $e\rond \theta = \sigma\rond e$,
where $\sigma(\alpha,\beta) = (\alpha,\beta)$. If we put $H = \ell^2(\Vk)$, $K
= \ell^2(\Ek)$ and the graph is locally finite, the graph structure yields
bounded operators $\Theta : K \to K$, $\un_w \mapsto \un_{\theta(w)}$ and $E : K
\to H\tens H$, $\un_w \mapsto \un_{e(w)}$. It is useful to consider
non-oriented, or ``geometrical'', edges obtained by identifying each edge with
the reversed one. At the hilbertian level we will use the subspace of
antisymmetric vectors $K_g = \Ker(\Theta + \id)$.

If moreover the graph is connected and endowed with an origin we have a natural
length function on vertices given by the distance to the origin and we can
consider ascending edges, i.e. edges which start at some length $n$ and end at
length $n+1$. At the hilbertian level we will use the partition of unity $\id_H
= \sum p_n$ given by the characteristic functions of the corresponding spheres,
and the orthogonal projection $p_+ : K \to K_+$ onto the subspace of functions
supported on ascending edges.

Denote finally by $\Hh$, $\Kk$, and $\Kk_g = \Kk \cap K_g$ the natural dense
subspaces of functions with finite support, and define $\epsilon : \Hh \to \CC$
by putting $\epsilon(\un_v) = 1$ for all $v\in\Vk$. Then one can introduce
target and source maps $\O = (\id\tens\epsilon) \rond E$, $\B =
(\epsilon\tens\id)\rond E$ which are bounded in the locally finite case. It is
an exercise to show that the graph under consideration is a tree \iff the
restriction $\B : \Kk_g \to \Ker\epsilon \subset \Hh$ is a bijection. The
corresponding issue in the quantum case will prove to be quite intricate.

\bigskip

Consider now a discrete group $\Gamma$ and fix a finite subset
$\Dd\subset\Gamma$ such that $\Dd^{-1} = \Dd$ and $e\notin\Dd$. The Cayley graph
associated with $(\Gamma,\Dd)$ can be described in two equivalent ways:
\begin{enumerate}
\item we take $\Vk = \Gamma$ and $\Ek = \{(g,h) \in \Gamma^2 ~|~ h\in g\Dd\}$,
  with the evident maps $\theta : (g,h) \mapsto (h,g)$ and $e : (g,h) \mapsto
  (g,h)$ ;
\item we take $\Vk = \Gamma$ and $\Ek = \Gamma\times \Dd$, with the reversing
  map $\theta : (g,h) \mapsto (gh, h^{-1})$ and the endpoints map $e : (g,h)
  \mapsto (g, gh)$.
\end{enumerate}
The first presentation leads to the notion of classical Cayley graph for
discrete quantum groups, whereas the second presentation, or rather the Hilbert
spaces and operators associated with it as above, lead to the notion of quantum
Cayley graph for discrete quantum groups \cite{Vergnioux_Cayley}.

\section{Generalities on path cocycles}
\label{sec_cocycles}

\subsection{Trivial and bounded cocycles}

Let $\Gamma$ be a discrete group, and $\Ss = \CC\Gamma$ the group algebra,
endowed with its natural structure of $*$-Hopf algebra. Extending maps from
$\Gamma$ to $\CC\Gamma$ by linearity, unitary representations $\pi : \Gamma \to
B(K)$ correspond to non-degenerate $*$-representations of $\Ss$, and
$\pi$-cocycles $c : \Gamma \to K$ yield derivations from $\Ss$ to the
$\Ss$-bimodule ${}_\pi K_\epsilon$, where $\epsilon$ is the trivial character on
$\Gamma$.

For the quantum case we make the following evident Definition:

\begin{df}
  Let $(\Ss,\epsilon)$ be a unital algebra with fixed character
  $\epsilon\in\Ss^*$. Let $\pi : \Ss \to B(K)$ be a unital representation on a
  vector space $K$. A {\em $\pi$-cocycle} is a derivation from $\Ss$ to the
  bimodule ${}_\pi K_\epsilon$, i.e. a linear map $c : \Ss \to K$ such that
  \begin{equation*}
    \forall x, y\in\Ss~~ c(xy) = \pi(x)c(y) + c(x)\epsilon(y).
  \end{equation*}
  We say that $c$ is {\em trivial} if the corresponding derivation is inner, i.e.
  there exists $\xi \in K$ such that $c(x) = \pi(x)\xi - \xi\epsilon(x)$. We
  call $\xi$ a {\em fixed vector} relatively to $c$.
\end{df}

When $(\Ss, \epsilon)$ is associated with a discrete quantum group as in
Section~\ref{sec_notation}, cocycles can also be considered from the dual point
of view: to any unital $*$-representation $\pi : \Ss \to L(K)$ is associated a
corepresentation $X \in M(\hat S\tens K(K))$, and to any $\pi$-cocycle $c : \Ss
\to K$ is associated an ``unbounded multiplier'' of the Hilbert \Cst module
$\hat S\tens K$:
\begin{displaymath}
  C = (\id\tens c)(V) \in (\hat S\tens K)^\eta 
  := \ts\prod_\alpha (p_\alpha\hat S\tens K),
\end{displaymath}
where $\alpha$ runs over $\Irr\Cc$. Note that $C$ is best understood as an
unbounded operator from $H$ to $H\tens K$. 

It is immediate to rewrite the cocyle relation as follows --- an equality
between elements of the algebraic tensor product $p_\alpha\hat S\tens
p_\beta\hat S\tens K$ after multiplying by any $p_\alpha\tens p_\beta\tens\id$
on the left:
\begin{displaymath}
 (\hat\delta\tens\id)(C) = V_{12}^*C_{23} V_{12} = X_{13} C_{23} + C_{13}.
\end{displaymath}
Similarly, $c$ is trivial with fixed vector $\xi$ \iff $C = X(1\tens\xi) -
1\tens\xi$. 

In this paper we will mainly work with $c : \Ss\to K$, but $C$ can be useful
when it comes to boundedness. We say that the cocycle $c$ is bounded if $C =
(\id\tens c)(V)$ is bounded in the following equivalent meanings : the family
$((p_\alpha\tens c)(V))_\alpha$ is bounded with respect to the \Cst hilbertian
norms ; the operator $(\id\tens c)(V) : \Hh \to H\tens K$ is bounded ; left
multiplication by $(\id\tens c)(V)$ defines an element of the $M(\hat
S)$-Hilbert module $M(\hat S\tens K)$ \cite{BaajSkand_KKS}.

\bigskip

It is a classical feature of cocycles with values in unitary representations
that triviality is equivalent to boundedness. We prove now a quantum
generalization of this result, which will not be used in the rest of the article
but we find is of independent interest.

The classical proof relies on the center construction for bounded subsets of
Hilbert spaces. In the quantum case, it is more convenient to rephrase it in
terms of metric projections. Recall that a metric projection of an element
$\zeta$ of a Banach space $E$ onto the subspace $F$ is an element $\xi\in F$
such that $\|\zeta-\xi\| = d(\zeta, F)$. The subspace $F$ is called Chebyshev if
all elements of $E$ have a unique metric projection onto $F$. All subspaces of
Hilbert spaces are Chebyshev, and in the case of Hilbert \Cst modules we have
the elementary:

\begin{lem} \label{lem_approx} Let $A$ be a unital \Cst algebra and $E$ a
  Hilbert $A$-module. Let $F \subset E$ be a closed subspace. Let $B$ be a \Cst
  algebra.
  \begin{enumerate}
  \item Assume $F$ is reflexive and $(\zeta|\zeta)$ is invertible in $A$ for
    every $\zeta \in F\setminus\{0\}$. Then $F$ is a Chebyshev subspace of $E$.
  \item Assume that $\xi\in E$ has a unique metric projection $\zeta\in F$ and
    consider the Hilbert $M(B\tens A)$-module $M(B\tens E)$. Then $1\tens \zeta$
    is the unique metric projection of $1\tens\xi$ onto $M(B\tens F)$.
  \end{enumerate}
\end{lem}

\begin{dem}
  1. Existence of metric projections in $F$ follows from reflexivity by weak
  compactness of closed balls in $F$ \cite[Cor.~2.1]{Singer_Approx}. Assume
  $\xi$, $ \xi' \in F$ are metric projections of $\zeta$. Then we put $d =
  d(\zeta, F)$, $\eta = (\xi+\xi')/2 \in F$ and we have, by the parallelogram
  identity for Hilbert modules :
  \begin{displaymath}
    (\zeta-\eta|\zeta-\eta) = \frac 12 (\zeta-\xi|\zeta-\xi) + 
    \frac 12 (\zeta-\xi'|\zeta-\xi') - \frac 14 (\xi-\xi'|\xi-\xi').
  \end{displaymath}
  Now we have $(\zeta-\xi|\zeta-\xi)$, $(\zeta-\xi'|\zeta-\xi') \leq d^2 1_A$
  and if $\xi \neq \xi'$ we get $\epsilon > 0$ such that $(\xi-\xi'|\xi-\xi')
  \geq \epsilon 1_A$, by hypothesis on $F$. Hence we get $\|\zeta-\eta\|^2 \leq
  d^2 - \epsilon /4 < d^2$, a contradiction.

  2. We clearly have $\|1\tens\zeta - 1\tens\xi\| = \|\zeta-\xi\| =:d$. Take
  $\eta\in M(B\tens F)$ such that $\|\eta - 1\tens\xi\| \leq d$. For any state
  $\varphi$ of $B$ there is a well-defined contraction $\varphi\tens\id :
  M(B\tens E) \to E$.  In particular $\|(\varphi\tens\id)(\eta)- \xi\| \leq d$,
  so the hypothesis implies $(\varphi\tens\id)(\eta) = \zeta$ for all
  $\varphi$. This clearly implies $\eta = 1\tens\zeta$.
\end{dem}

\begin{prp} \label{prp_bounded_trivial}
  Let $\pi : S \to L(K)$ be a $*$-representation of the full \Cst algebra of a
  discrete quantum group, and $c : \Ss \to K$ a $\pi$-cocycle. If $(\id\tens
  c)(V) \in M(\hat S\tens K)$ then $c$ is trivial.
\end{prp}

\begin{dem}
  Put $C = (\id\tens c)(V) \in M(\hat S\tens K)$ and let $X \in L(H\tens K)$ be
  the corepresentation corresponding to $\pi$. We apply Lemma~\ref{lem_approx}
  to the subspace $F = 1\tens K$ of $E = M(\hat S\tens K)$. Reflexivity is
  immediate since $F$ is a Hilbert space, and the second condition of the Lemma
  is trivially satisfied since $(1\tens \zeta| 1\tens \zeta) = \|\zeta\|^2_K\,
  1_{M(\hat S)}$. Denote by $1\tens \xi$ the metric projection of $C$ onto $1\tens
  K$ and put $d = \|1\tens\xi - C\|$.

  By Lemma~\ref{lem_approx} the vector $1\tens 1\tens\xi$ is the metric
  projection of $1\tens C$ onto the subspace $M(\hat S\tens 1\tens K)$. Since
  $X_{13}$ is a unitary which stabilizes $M(\hat S \tens 1\tens K)$ and $C_{13}$
  belongs to $M(\hat S\tens 1\tens K)$, we can deduce that $V_{12}^* C_{23}
  V_{12} = X_{13} C_{23} + C_{13}$ lies at distance $d$ from $M(\hat S \tens
  1\tens K)$ with unique metric projection $X_{13}(1\tens 1\tens\xi) +
  C_{13}$. But $1\tens 1\tens\xi$ is evidently a vector in $M(\hat S \tens
  1\tens K)$ lying at distance $d$ from $V_{12}^* C_{23} V_{12}$. Hence we have
  that $X(1\tens \xi) + C = 1\tens\xi$, and the cocycle is trivial.
\end{dem}

\subsection{Path cocycles}

We investigate in this paper a geometrical method to construct particular
cocycles on $\Ss$, the so-called path cocycles. Before introducing this notion
we need to recall the definitions of the classical and quantum Cayley graphs
associated with a discrete quantum group.

\begin{rque}{Reminder} \label{rem_cayley} As explained in
  Section~\ref{sec_notation}, let $\Cc$ be the category of finite-dimen\-sional
  representations of $\hat S$, and let $\Irr\Cc$ be a system of representants of
  all irreducible objects. We fix a finite subset $\Dd \subset \Irr\Cc$ such
  that $\bar \Dd = \Dd$ and $1_\Cc\notin\Dd$. Recall from
  \cite[Def.~3.1]{Vergnioux_Cayley} that the classical Cayley graph associated
  with $(\Ss,\delta,\Dd)$ is given by:
  \begin{itemize}
  \item the set of vertices $\Vk = \Irr \Cc$ and the set of edges 
    \begin{displaymath}
      \Ek = \{(\alpha,\beta) \in (\Irr\Cc)^2 ~|~ \exists \gamma\in\Dd~ 
      \beta\subset \alpha\tens\gamma\},
    \end{displaymath}
  \item the canonical reversing map $\sigma : \Ek \to \Ek, (\alpha, \beta)
    \mapsto (\beta, \alpha)$ and endpoints map $i_{\mathrm{can}} : \Ek \to
    \Vk\times \Vk$.
  \end{itemize}
  We endow this graph with the root $1_\Cc$ given by the trivial
  corepresentation. The elements of $\Dd$ are called directions of the Cayley
  graph. In the classical case of a discrete group $\Gamma$ we recover the first
  description of the usual Cayley graph given in Section~\ref{sec_notation}.

  On the other hand, let $(H,V,U)$ be the Kac triple associated with
  $(S,\delta)$ and denote by $\Sigma \in B(H\tens H)$ the flip operator.  We
  introduce the central projection $p_1 \in \hat S$ corresponding to $\Dd
  \subset \Irr \Cc$ as explained in Section~\ref{sec_notation}. Recall from
  \cite[Def.~3.1]{Vergnioux_Cayley} that the hilbertian quantum Cayley graph
  associated with $(\Ss,\delta,p_1)$ is given by:
  \begin{itemize}
  \item the $\ell^2$-space of vertices $H$ and the $\ell^2$-space of edges $K =
    H\tens p_1 H$,
  \item the reversing operator $\Theta = \Sigma (1\tens U) V (U\tens U) \Sigma
    \in B(K)$ and the endpoints operator $\E = V_{|K} : K \to H\tens H$.
  \end{itemize}
  In the classical case of a discrete group $\Gamma$ one recovers the hilbertian
  objects associated with the usual Cayley graph, as introduced in
  Section~\ref{sec_notation}. A major feature of the quantum case is that
  $\Theta$ need not to be involutive, even in a ``deformed'' meaning.

  In fact we will rather use the source and target maps $\O = (\id\tens\epsilon)
  V : K \to H$ and $\B = (\epsilon\tens\id) V : K \to H$, which are bounded if
  $\Dd$ is finite.  We denote by $K_g = \Ker (\Theta+\id)$ the space of
  geometrical, or antisymmetric, edges and by $p_g \in B(K)$ the orthogonal
  projection onto $K_g$. The spaces $H$, $K$, $K_g$ are naturally endowed with
  representations of $S_\red$ which are intertwined by $\O$, $\B$ and
  $\Theta$.
\end{rque}

Recall that we have a natural dense subspace $\Hh = \Lambda_h(\Ss)$ of
$H$, and define similarly $\Kk = \Hh\tens p_1 H$ at the level of edges. The
situation is more complicated for geometrical edges and we make the following
definitions, which will be discussed in Remark~\ref{rk_comment_dense}. We call
$\Kk_g$ the intersection of $K_g$ with $\Kk$, and we call $\Kk'_g$ the
orthogonal projection of $\Kk$ onto $K_g$.

\begin{df} \label{df_path_cocycle}
  The {\em trivial $\lambda$-cocycle} of $\Ss$ is the cocycle with fixed vector
  $\xi_0$, i.e. $c_0 : x \mapsto \Lambda_h(x) - \epsilon(x)\Lambda_h(1)$. A
  cocycle $c_g : \Ss \to K$ is called a {\em path cocycle} if we have $\B\rond
  c_g = c_0$ and $c(\Ss) \subset \Kk'_g$.
\end{df}

The motivating example is as follows: for any element $g$ of the free group
$F_n$, consider the unique path from $e$ to $g$ in the Cayley graph of $F_n$,
and let $c_g(g) \in K_g$ denote the sum of (characteristic functions of)
antisymmetric edges along this oriented path. Then it is easy to check that
$c_g$ extends by linearity to a path cocycle on $\CC F_n$. The evident fact that
this cocycle is proper (as a function on $F_n$) establishes the a-T-menability
of $F_n$.

The problem of existence of path cocyles is a problem about the injectivity and
surjectivity of appropriate restrictions of $\B : K_g \to H$. As indicated in
Section~\ref{sec_notation}, it is thus related to the question of whether the
quantum Cayley graph under consideration can be considered a tree. Notice that
we considered in \cite{Vergnioux_Cayley} another approach to this question in
terms of the ascending orientation.

We start by attacking the surjectivity side of the issue with the simple
Lemma~\ref{lem_B-O} below. Recall that we consider the counit $\epsilon$ as a
linear map $\epsilon : \Hh \to \CC$, and its kernel $\Ker \epsilon$ as a
subspace of $\Hh$. We clearly have $\Ker \epsilon = \{\xi - \epsilon(\xi)\xi_0
~|~ \xi\in\Hh\} = c_0(\Ss)$. Hence for our purposes ``surjectivity'' of $\B :
K_g \to H$ corresponds to the requirement that $\B(\Kk'_g)$ contains
$\Ker\epsilon$.

\begin{rque}{Reminder} \label{rem_ascending} 
  To state and prove the Lemma we need to recall the notion of ascending
  orientation $K_{\ss++}\subset K$ for quantum Cayley graphs.

  We begin with spheres. The distance to the origin in the classical Cayley
  graph yields spheres of radius $n$ which are subsets of $\Irr\Cc$. To these
  subsets are associated central projections $p_n \in \hat S \subset
  B(H)$. These projections play the role of characteristic functions of spheres
  for the quantum Cayley graph. If the classical Cayley graph is connected, they
  form a partition of the unit of $B(H)$.

  Notice that the notation $p_n$ is consistent with $p_0$, $p_1$ introduced
  earlier. When the space acted upon is clear from the context, we will moreover
  denote $p_n = p_n\tens\id \in B(K)$, which projects onto ``edges starting at
  distance $n$ from the origin''. We will also use $p_{<n} = \sum_{k<n} p_k$,
  $p_{\geq n} = \sum_{k\geq n} p_k$ etc.

  Next we put $p_{\ss\jok+} = \sum (p_n\tens p_1)\hat\delta(p_{n+1})$ and
  $p_{\ss+\jok} = \sum (p_n\tens p_1)\hat\delta'(p_{n+1})$, where
  $\hat\delta'(a) = (U\tens U)\Sigma\delta(a)\Sigma(U\tens U) \in B(H\tens
  H)$. The projection $p_{\ss\jok+}$ lives in $\hat S\tens\hat S \subset
  B(H\tens H)$, and $p_{\ss+\jok}$ is the corresponding projection in $U\hat
  SU\tens U\hat SU$, acting ``from the right'' on $H$. The ascending projection
  of the quantum Cayley graph is by definition $p_{\ss++} =
  p_{\ss\jok+}p_{\ss+\jok} \in B(K)$. Its image is denoted $K_{\ss++} =
  p_{\ss++}K$, and we put $\Kk_{\ss++} = \Kk \cap K_{\ss++} = p_{\ss++} \Kk$.

  In the classical case it is easy to check that $p_{\ss+\jok} = p_{\ss\jok+} =
  p_+$ is indeed the orthogonal projection onto the subspace of functions
  supported on ascending edges. Similarly $p_n$ is the orthogonal projection
  onto the subspace of functions supported on the sphere of radius $n$.
\end{rque}

\begin{lem} \label{lem_B-O} Assume that the classical Cayley graph associated
  with $(\Ss, \Dd)$ is connected, i.e. $\Dd$ generates $\Cc$. We have
  then the following results for the quantum Cayley graph.
  \begin{enumerate}
  \item We have $(\B-\O)(K_g^\bot) = 0$.
  \item We have $(\B-\O) (\Kk_{\ss++}) = \Ker \epsilon$.
  \end{enumerate}
  As a result we have $\B(\Kk_g) \subset \Ker \epsilon = \B(\Kk'_g) \subset
  \B(K_g)$.
\end{lem}

\begin{dem}
  1. Using the identities $\O = \B\Theta = \B\Theta^*$
  \cite[Prop.~3.6]{Vergnioux_Cayley} we have:
  \begin{displaymath}
    (\Theta+\id)(\B-\O)^* = (\Theta+\id)(\id-\Theta^*)\B^* =
    (\Theta-\Theta^*)\B^* = 0.   
  \end{displaymath}
  Hence $\Img(\B-\O)^* \subset K_g$ and $K_g^\bot \subset \Ker(\B-\O)$.

  Let us explain also the last statement. We have $\B
  (\Lambda_h\tens\Lambda_h)(x\tens y) = \Lambda_h(xy)$ and
  $\O(\Lambda_h\tens\Lambda_h)(x\tens y) = \Lambda_h(x \epsilon(y))$ for $x$,
  $y\in\Ss$ \cite[Prop.~3.6]{Vergnioux_Cayley}. Since $\epsilon(xy) =
  \epsilon(x\epsilon(y))$, this shows that $(\B-\O)(\Kk)$, hence $\B(\Kk_g) =
  (\B-\O)(\Kk_g)$ and $\B(\Kk'_g)$, are contained in $\Ker\epsilon$. On the
  other hand 1. and 2. clearly imply that $\B(\Kk'_g) = \Ker \epsilon$, by
  definition of $\Kk'_g$.

  2. The previous paragraph already shows that $(\B-\O)(\Kk_{\ss++}) \subset
  \Ker \epsilon$. Hence it remains to prove that $\xi - \epsilon(\xi)\xi_0 \in
  (\B-\O) (\Kk_{\ss ++})$ for all $\xi \in p_nH$ and all $n$, and we will
  proceed by induction over $n$.

  We will need to know that $\B (p_n K_{\ss++}) = p_{n+1} H$. The proof is as
  follows: for every $\gamma \in \Irr\Cc$ of length $n+1$ there exists by
  definition of the length $\alpha$, $\beta \in \Irr\Cc$ of respective lengths
  $n$, $1$ such that $\gamma \subset \alpha\tens\beta$. Then the relevant
  arguments of \cite[Prop.~4.7 and Rem.~4.8]{Vergnioux_Cayley} apply in our more
  general setting and show that $\B : \hat\delta(p_\gamma) (p_\alpha\tens
  p_\beta) K_{\ss++} \to p_\gamma H$ is an explicit, non-zero multiple of a
  (surjective) isometry.

  We have now $\B(p_n K_{\ss ++}) = p_{n+1}H$ and $\O(p_nK_{\ss++}) \subset
  p_nH$. For a fixed $\xi\in p_n H$, this yields $\zeta\in p_{n-1}K_{\ss++}$
  such that $\B\zeta = \xi$, and $\xi' = \O\zeta\in p_{n-1}H$. Since $\B-\O$
  ranges in $\Ker \epsilon$, we have $\epsilon(\xi) = \epsilon(\xi')$. Moreover,
  by induction over $n$ one can find $\zeta' \in \Kk_{\ss++}$ such that
  $(\B-\O)(\zeta') = \xi' - \epsilon(\xi')\xi_0$. We have then
  \begin{equation*}
    \xi-\epsilon(\xi)\xi_0 = (\B-\O)(\zeta) + \xi' - \epsilon(\xi')\xi_0 = 
    (\B-\O) (\zeta+\zeta').
  \end{equation*}
  For $n=0$ the result is obvious since $p_0 H = \CC \xi_0$ and $\xi_0 -
  \epsilon(\xi_0)\xi_0 = 0$.
\end{dem}

\begin{crl} \label{crl_path_prop} If $c_g$ is a path cocyle, then $c_g(\Ss)$ is
  an $\Ss$-stable subspace of $K_g$ on which $\B$ is injective. Conversely, let
  $K'_g$ be an $\Ss$-stable subspace of $K_g$ on which $\B$ is injective, and
  assume moreover that $\Kk'_g \subset K'_g$. Then there exists a unique path
  cocyle $c_g$ with values in $K'_g$, and $c_g(\Ss) = \Kk'_g$.
\end{crl}

\begin{dem}
  The cocyle relation implies trivially the $\Ss$-stability of $c_g(\Ss)$.  Assume
  $\B(\zeta) = 0$ with $\zeta = c_g(x)$, $x\in\Ss$. Then we have $0 = \B(\zeta) =
  c_0(x)$ by definition of a path cocyle. But clearly $c_0(x) = 0$ \iff $x \in
  \CC 1$, and since $1$ is mapped to $0$ by any cocycle (take $x=y=1$ in the
  cocycle relation) we conclude that $\zeta = 0$.

  If $K'_g$ contains $\Kk'_g$, then $\B(K'_g)$ contains $\Ker \epsilon =
  c_0(\Ss)$ by Lemma~\ref{lem_B-O}. If moreover $\B$ is injective on $K'_g$ we
  obtain a unique map $c_g : \Ss \to K'_g$ such that $\B c_g = c_0$, and its
  image equals $\Kk'_g$. By uniqueness, $\Ss$-stability, and since $\B$
  intertwines the actions of $\Ss$, the cocycle relation for $c_0$ implies the
  cocycle relation for $c_g$.
\end{dem}

\begin{rque}{Remark} {\sc Comment on $\Kk'_g$.} \\ \label{rk_comment_dense} In
  the classical case, both subspaces $\Kk_g$ and $\Kk'_g$ coincide with the
  space of functions with finite support on antisymmetric edges. In the quantum
  case, the inclusion $\Kk_g \subset \Kk'_g$ is strict in general : we will
  e.g. argue in Section~\ref{sec_endpoints} that $\Kk'_g$ is not included in
  $\Kk$, whereas this is by definition the case of $\Kk_g$.

  Note that $\Kk_g$ and $\Kk'_g$ are both stable under the action of the dense
  subalgebra $\Ss \subset S$, because this is the case of $\Kk$ and $K_g$. The
  subspace $\Kk'_g$ is obviously dense in $K_g$, whereas this is by no mean
  clear for $\Kk_g$ --- we give however a proof of this density at the end of
  Theorem~\ref{thm_orient_ana}.

  The choice of $\Kk'_g$ instead of $\Kk_g$ in Definition~\ref{df_path_cocycle}
  is motivated by the results of Lemma~\ref{lem_B-O} and
  Corollary~\ref{crl_path_prop}. Notice that these ``surjectivity'' results
  still hold if one considers $\Kk''_g = p_g(\Kk_{\ss++})$ instead of $\Kk'_g$,
  and this smaller subspace is a priori more suited for injectivity
  results. However the $\Ss$-stability and density of $\Kk''_g$ are not clear
  anymore. In fact one can prove, using the results in this article, that
  $\Kk'_g = \Kk''_g$ for free products of universal quantum groups.
\end{rque}

\section{Complements on quantum Cayley trees}
\label{sec_endpoints}

In this section we give a new description of the closed subspace $p_{\ss++}K_g$
already studied in \cite{Vergnioux_Cayley} and we prove a density result for
$\Kk_g \subset K_g$. The main new tools are a ``rotation operator'' $\Phi :
K_{\ss--} \to K_{\ss++}$, the associated shift $r' : K_{\ss++} \to K_{\ss++}$,
and a ``quasi-classical'' subspace $Q_0K$.

\begin{rque}{Reminder} \label{rem_alg} We assume now that the discrete quantum
  group under consideration is a free product of orthogonal and unitary
  universal discrete quantum groups, and that $\Dd$ is the collection of the
  corresponding fundamental corepresentations $(u_{ij})_{ij}$ --- this is the
  setting of \cite{Vergnioux_Cayley}.

  In this case the classical Cayley graph is a tree, and for each edge $w =
  (\alpha,\beta)$ one can define the direction of $w$ as the unique $\gamma \in
  \Dd$ such that $\beta\subset \alpha\tens\gamma$, or equivalently,
  $(p_\alpha\tens p_\gamma) \hat\delta(p_\beta) \neq 0$. Notice that there can
  be (at most) two edges starting from $\alpha$ in direction $\gamma$, one of
  them ascending and the other one descending.

  If one puts $p_{\ss\jok-} = \sum (p_n\tens p_1)\hat\delta(p_{n-1})$ and
  $p_{\ss-\jok} = \sum (p_n\tens p_1)\hat\delta'(p_{n-1})$, then
  $\{p_{\ss\jok+},$ $p_{\ss\jok-}\}$ and $\{p_{\ss+\jok}, p_{\ss-\jok}\}$ are
  commuting partitions of unity, which also commute to the partition $\{p_n\}$.
  We have the following related pieces of notation: $p_{\ss+-} =
  p_{\ss+\jok}p_{\ss\jok-}$, $K_{\ss+-} = p_{\ss+-}K$, $\Kk_{\ss+-} = \Kk\cap
  K_{\ss+-}$ and similarly for $p_{\ss-+}$ and $p_{\ss--}$.

  The behaviors of $\B$ and $\Theta$ with respect to these partitions of unity
  yield ``computation rules'' in quantum Cayley trees, cf.
  \cite[Prop.~4.3]{Vergnioux_Cayley}, \cite[Prop.~5.1]{Vergnioux_Cayley} and the
  identity $\B \Theta = \O = \B \Theta^*$ which results from
  \cite[Prop.~3.6]{Vergnioux_Cayley}. The interplay between orientation and
  reversing yields two important operators: $r = - p_{\ss+-}\Theta p_{\ss+-}$
  and $s = p_{\ss+-}\Theta p_{\ss++}$.
\end{rque}

\subsection{New tools}

In the following Lemma we introduce a right shift $r'$ on $K_{\ss++}$ which will
be used in the statement and the proof of Theorem~\ref{thm_orient_ana}. This is
done via the rotation operator $\Phi$ which is of independent interest. In the
case of a classical tree, this operator corresponds to ``rotating'' descending
edges around their target vertex, yielding the unique ascending edge with the
same target.

\begin{lem} \label{lem_rot} Consider the quantum Cayley graph of a free product
  of universal discrete quantum groups.
  \begin{enumerate}
  \item There exists a unique bounded operator $\Phi : K_{\ss--} \to K_{\ss++}$ such
    that \newline $(1-p_0) \B p_{\ss--} = \B p_{\ss++} \Phi$. We have $\Phi(p_n
    K_{\ss--}) \subset p_{n-2}K_{\ss++}$ for $n\geq 2$ and $\Phi(p_1 K_{\ss--})=\{0\}$.
  \item The operator $p_{\ss++}\Theta p_{\ss--}\Phi^*$ is a right shift $r'$ on
    $K_{\ss++}$ such that $sr'=rs$. Moreover we have $\Phi p_{\ss--}\Theta
    p_{\ss+-} = s^*$.
  \end{enumerate}
\end{lem}

\begin{dem} 1. We know that ${\B}_{|K_{\ss++}} : K_{\ss++} \to H$ is injective
  and that its image equals $(1-p_0)H$ \cite[Prop.~4.7]{Vergnioux_Cayley}, and
  we denote $\uB = \B : K_{\ss++} \to (1-q_0)H$ the corresponding invertible
  restriction. Moreover we have $(1-p_0)\B p_{\ss--} = \B
  p_{--}(1-p_1)$. Therefore it suffices to put $\Phi = \uB^{-1}(\B
  p_{\ss--})(1-p_1)$. The remaining statements hold because $\B p_n K_{\ss--} =
  p_{n-1} H$.

  2. Using the computation rules in quantum Cayley trees, we can write
  \begin{align} \nonumber
    \B p_{\ss++} (rs)^* &= -\B p_{\ss++}\Theta^*p_{\ss+-}\Theta^* p_{\ss+-}
    = -\B p_{\ss+\jok} \Theta^*p_{\ss+-}\Theta^* p_{\ss+-} \\ \nonumber
    &= -\B \Theta^*p_{\ss+-}\Theta^* p_{\ss+-}
    = -\B \Theta p_{\ss+-}\Theta^* p_{\ss+-} \\ \label{eq_shifts_1}
    &= -\B p_{\ss--}\Theta p_{\ss+-}\Theta^* p_{\ss+-}.
  \end{align}
  Similarly, using the definition of $\Phi$:
  \begin{align} \nonumber
    \B p_{\ss++} (sr')^* 
    &= \B p_{\ss++} \Phi p_{\ss--} \Theta^* p_{\ss++} \Theta^* p_{\ss+-} \\
    &= (1-p_0)\B p_{\ss--} \Theta p_{\ss++}\Theta^* p_{\ss+-}
    = \B p_{\ss--} \Theta p_{\ss++}\Theta^* p_{\ss+-}. \label{eq_shifts_2}
  \end{align}
  Notice indeed that $p_0 \B p_{\ss--} \Theta p_{\ss++}\Theta^* p_{\ss+-}$
  vanishes, essentially because ``all edges starting from the
  origin are ascending'', i.e. $p_0 p_{\ss++} = p_0$:
  \begin{align*}
    p_0 \B p_{\ss--} \Theta p_{\ss++}\Theta^* p_{\ss+-} &=
    \B p_{\ss--} p_1 \Theta p_{\ss++}\Theta^* p_{\ss+-} =
    \B p_{\ss--} \Theta p_0p_{\ss++}\Theta^* p_{\ss+-} \\ &=
    \B p_{\ss--} p_1 \Theta \Theta^* p_{\ss+-} = \B p_1 p_{\ss--} p_{\ss+-} = 0.
  \end{align*}
  Now we subtract~\eqref{eq_shifts_1} from~\eqref{eq_shifts_2}:
  \begin{align*}
    \B p_{\ss++} (sr'-rs)^* &= \B p_{\ss--}\Theta p_{\ss+\jok}\Theta^* p_{\ss+-} \\
    &= \B p_{\ss--}\Theta\Theta^*p_{\ss+-} = \B p_{\ss--} p_{\ss+-} = 0.
  \end{align*}
  Since $\B$ is injective on $K_{\ss++}$, this yields $sr'=rs$. The identity
  $\Phi p_{\ss--}\Theta p_{\ss+-} = s^*$ follows similarly from
  \begin{align*}
    \B p_{\ss++} \Phi p_{\ss--}\Theta p_{\ss+-} &= \B p_{\ss--} \Theta p_{\ss+-}
    = \B \Theta p_{\ss+-} = \B \Theta^* p_{\ss+-} \\
    &= \B p_{\ss++} \Theta^* p_{\ss+-} = \B p_{\ss++} s^*.
  \end{align*}\par
\end{dem}

We introduce now quasi-classical subspaces $Q_0K$, $Q_0H$, or rather the
orthogonal projection $Q_0$ onto them. The name ``quasi-classical'' refers to
the fact that we will have $\Theta^2 Q_0 = Q_0$ and $(p_{\ss+-}+ p_{\ss-+}) Q_0
= 0$, as in classical trees.

In the case of a single copy of an orthogonal group $A_o(Q)$, the projection
$Q_0$ coincides in fact with $q_0$ introduced in
\cite[p.~125]{Vergnioux_Cayley}, but we will give here an independent
definition. In general, $Q_0 K$ contains $q_0K$ strictly and we call $q_0 K$ and
$q_0 H$ the ``classical subspaces'' of the quantum Cayley tree, because they are
really the $\ell^2$-spaces of the classical Cayley tree --- see the beginning of
Section~\ref{sec_cayley} for more details.

Note that $Q_0H$, $Q_0K$ do not really define a (quantum) subtree of $(H,K)$
because $\B Q_0 K$ is not included in $Q_0 H$ in general. Nevertheless, we prove
below that $\B Q_0 K_{\ss++} \subset Q_0 H$, and also $\B (1-Q_0) K \subset
(1-Q_0)H$, i.e. the ``purely quantum'' subspaces $(1-Q_0)H$, $(1-Q_0)K$ form a
quantum subtree of $(K,H)$.

\begin{lem}\label{lem_quasiclass} 
  Consider the quantum Cayley graph of a free product of universal discrete
  quantum groups.
  \begin{enumerate}
  \item Denote by $Q_0 \in B(K)$ the orthogonal projection onto
    \begin{displaymath}
      Q_0K = \{\zeta \in K_{\ss++} ~|~ \Theta\zeta\in K_{\ss--}\} \oplus
      \{\zeta \in K_{\ss--} ~|~ \Theta\zeta\in K_{\ss++}\}.
    \end{displaymath}
    Then we have $[Q_0, p_n] = [Q_0, p_{\ss+\jok}] = [Q_0, p_{\ss\jok+}] = [Q_0,\Theta]
    = 0$. \\ Moreover $\Theta^2Q_0 = Q_0$ and $Q_0 K_{\ss++} = p_{\ss++} \Ker s$.
  \item Denote by $Q_0 \in B(H)$ the orthogonal projection onto
    \begin{displaymath}
      Q_0H = \{\zeta\in H ~|~ \O^*\zeta \in K_{\ss++}\oplus K_{\ss--}\}.
    \end{displaymath}
    Then $[Q_0, p_n] = 0$ and we have $Q_0\B(1-Q_0) = (1-Q_0)\B Q_0 p_{\ss++} =
    0$.
\end{enumerate}
\end{lem}

\begin{dem}
  1. It is clear from the definition that $Q_0$ commutes with $p_{\ss\jok+}$ and
  $p_{\ss+\jok}$. Since $\Theta p_{\ss++} p_n = p_{n+1}\Theta p_{\ss++}$, and
  similarly with $p_{\ss--}$, it is also easy to check that $Q_0$ commutes with
  $p_n$.

  Since $\Theta p_{\ss++} = p_{\ss\jok-}\Theta p_{\ss++}$ and $\Theta p_{\ss--}
  = p_{\ss\jok+}\Theta p_{\ss--}$, the definition of $Q_0K$ can be rewritten
  \begin{displaymath}
    Q_0K = p_{\ss++} \Ker(p_{\ss+-}\Theta p_{\ss++}) \oplus
    p_{\ss--} \Ker(p_{\ss-+}\Theta p_{\ss--}).
  \end{displaymath}
  In particular we see that $Q_0K_{\ss++} = \Ker s$. Let us notice that the
  first term in the direct sum above also equals $p_{\ss++}
  \Ker(p_{\ss-+}\Theta^*p_{\ss++})$. We can indeed use the operator $W$ from
  \cite[Lem.~5.2]{Vergnioux_Cayley} to write, for any $\zeta \in K_{\ss++}$:
  \begin{align*}
    p_{\ss+-}\Theta p_{\ss++}\zeta = 0 ~~\Longleftrightarrow~~
    W p_{\ss+-}\Theta p_{\ss++}\zeta = 0 ~~\Longleftrightarrow~~
    p_{\ss-+}\Theta^* p_{\ss++}\zeta = 0.
  \end{align*}
  Similarly, the second term in the direct sum equals $p_{\ss--}
  \Ker(p_{\ss+-}\Theta^*p_{\ss--})$. 

  Now we can show that $Q_0$ commutes with $\Theta$. Take indeed $\zeta \in Q_0
  K_{\ss++}$: then $\Theta\zeta$ lies in $K_{\ss--}$, and in fact in $Q_0
  K_{\ss--}$ because $p_{\ss+-}\Theta^*\Theta\zeta = p_{\ss+-}\zeta =
  0$. Similarly, we see that $\Theta Q_0 K_{\ss++}$, $\Theta^* Q_0 K_{\ss++}$
  and $\Theta^* Q_0 K_{\ss--}$ are all contained in $Q_0 K$, hence $[Q_0,
  \Theta] = 0$. Finally, since $\Theta$ and $\Theta^*$ send $Q_0 K_{\ss++}$ to
  $K_{\ss--}$ we can write, using \cite[Prop.~5.1]{Vergnioux_Cayley}:
  \begin{displaymath}
    \Theta^2 p_{\ss++} Q_0 = \Theta p_{\ss --} \Theta p_{\ss++} Q_0
    = \Theta p_{\ss --} \Theta^* p_{\ss++} Q_0 
    = \Theta \Theta^* p_{\ss++} Q_0 = p_{\ss++} Q_0. 
  \end{displaymath}
  Similarly $\Theta^2 p_{\ss--} Q_0 = p_{\ss--} Q_0$, hence $\Theta^2 Q_0 = Q_0$.

  2. We have $[Q_0, p_n] = 0$ because $\O(p_n\tens\id) = p_n \O$. By definition
  of $Q_0 K$ and remarks in the first part of the proof, we have 
  \begin{displaymath}
    (1-Q_0)K = K_{\ss+-}\oplus K_{\ss-+}\oplus \overline\Img\ p_{\ss++}\Theta p_{\ss-+} 
    \oplus \overline\Img\ p_{\ss--}\Theta p_{\ss+-}.
  \end{displaymath}
  It is known that $\B$ vanishes on $K_{\ss-+}$ and $K_{\ss+-}$
  \cite[Prop.~4.3]{Vergnioux_Cayley}. Consider then $\zeta = (p_{\ss++}\Theta
  p_{\ss-+} + p_{\ss--}\Theta p_{\ss+-}) \eta \in (1-Q_0)K$. We have 
  \begin{displaymath}
    \B \zeta = \B (\Theta p_{\ss-+} + \Theta p_{\ss+-}) \eta 
    = \O (p_{\ss-+} + p_{\ss+-}) \eta.
  \end{displaymath}
  Observing moreover that $Q_0 H = \Ker (p_{\ss+-}+p_{\ss-+})\O^* =
  (\O(K_{\ss+-} \oplus K_{\ss-+}))^\bot$, we obtain that $Q_0\B(1-Q_0) = 0$.

  The last property requires more care. Taking $\zeta \in Q_0 K_{\ss++}$, we
  have to show that $\O^*\B\zeta = \Theta \B^*\B\zeta$ lies in $K_{\ss++} \oplus
  K_{\ss--}$. Since $\B^*H \subset K_{\ss++}\oplus K_{\ss--}$, this is
  equivalent to the fact that $\B^*\uB\zeta$ lies in $Q_0 K$, where $\uB$
  denotes the invertible restriction $\uB = \B : K_{\ss++} \to (1-p_0)H$ as in
  the proof of Lemma~\ref{lem_rot}. First we want to replace $\uB$ with
  $(\uB^{-1})^*$. This is possible because $E_2$ is a multiple of an isometry on
  each subspace $(p_\alpha\tens p_\gamma) K_{\ss++}$ \cite[proof of
  Prop.~4.7]{Vergnioux_Cayley}, and $Q_0$ commutes not only with $(p_n \tens\id)
  p_{\ss++}$, but more precisely with each $(p_\alpha\tens p_\gamma) p_{\ss++}$.

  Let us prove now that $\B^*(\uB^{-1})^*\zeta$ lies in $Q_0K$. This is clearly
  the case for the $p_{\ss++}$ component since $\zeta\in Q_0K$ and
  \begin{displaymath}
    p_{\ss++}\B^*(\uB^{-1})^* = \uB^*(\uB^{-1})^* = \id_{K_{\ss++}}.
  \end{displaymath}
  For the $p_{\ss--}$ component, notice that $p_{\ss--}\B^*(\uB^{-1})^* =
  \Phi^*$. Hence we finally check that $\Phi^*\zeta$ lies in $Q_0 K$ using the
  last point of Lemma~\ref{lem_rot}: we have $p_{\ss+-}\Theta^* p_{\ss--}\Phi^*
  \zeta = s \zeta$, and since $\zeta \in Q_0K$ this vector vanishes.
\end{dem}

\subsection{Geometrical edges}

In the next Theorem we give a new description of the space $K_g$ basing on the
shift $r'$. This description is the key point for the proof of the injectivity
of $\B$ on $(1-Q_0)K_g$, see Proposition~\ref{prp_q_path}. Using the same
techniques, we also take the opportunity to prove the density of the subspace
$\Kk_g \subset K_g$.

Note that the non-involutivity of $\Theta$ in the quantum case, and more
precisely the fact that $p_{\ss+-}\Theta p_{\ss+-}$ acts as a shift with respect
to the distance to the origin, implies that the spectral projection $p_g$ does
not stabilize $\Kk$, and $\Kk'_g$ is not included in $\Kk$. In particular,
although $\Kk'_g = \Kk_g$ in the classical case, this is no longer true in the
quantum case, and Corollary~\ref{crl_path_prop} shows that $\Kk'_g$ is the
correct dense subspace for the study of path cocycles.

\begin{rque}{Reminder} \label{rem_Kg} Recall the following results from
  \cite{Vergnioux_Cayley}. The map $r = -p_{\ss+-}\Theta p_{\ss+-}$ acts as a
  right shift in the decomposition $K_{\ss+-} = \bigoplus p_kK_{\ss+-}$ and we
  denote by $K_\infty$ the inductive limit of the contractive system $(p_k
  K_{\ss+-}, r)$ of Hilbert spaces. Let $R_k : p_k K_{\ss+-} \to K_\infty$ and
  $R = \sum R_k : \Kk_{\ss+-} \to K_\infty$ be the canonical maps, and recall
  the notation $s = p_{\ss+-}\Theta p_{\ss++}$. Then the maps $R_k$ are
  injective, and the operator $Rs : K_{\ss++}\to K_\infty$ is a co-isometry
  \cite[Prop.~6.2 and Thm.~6.5]{Vergnioux_Cayley}.

  Besides it is not hard to prove that the restriction $p_{\ss++} : K_g \to
  K_{\ss++}$ is injective. When ``all directions have quantum dimensions
  different from~$2$'', we have the following identities \cite[Thms.~5.3
  and~6.5]{Vergnioux_Cayley}, which show in particular that $p_{\ss++}K_g$ is
  closed :
  \begin{equation} \label{eq_JV_image} 
    p_{\ss++}K_g = \{\zeta\in K_{\ss++} ~|~ \exists \eta\in K_{\ss+-}~ (1-r)\eta =
    s\zeta\} = \Ker Rs.
  \end{equation}
  Moreover if $\zeta = p_{\ss++}\xi$ and $\eta = - p_{\ss+-}\xi$ with $\xi \in
  K_g$, then $s\zeta = (1-r)\eta$ and $\zeta$, $\eta$ are related by the
  identities $R_kp_k\eta = Rs p_{< k}\zeta = -Rs p_{\geq k}\zeta$. Recall that in
  the classical case $K_{\ss+-} = K_\infty = 0$ so that $p_{\ss ++}K_g =
  K_{\ss++}$, as expected.

  Let us discuss rapidly the hypothesis that ``all directions have quantum
  dimensions different from $2$''. For a free products of quantum groups
  $A_o(Q)$, with $Q\bar Q \in \CC I_n$, and $A_u(Q)$, there is one direction
  $\gamma = \bar\gamma$ in the classical Cayley tree for each factor $A_o(Q)$,
  and two directions $\gamma$, $\bar\gamma$ for each factor $A_u(Q)$. These
  ``fundamental corepresentations'' have quantum dimension strictly greater than
  $2$, except in the following cases and isomorphic ones: $\dim_q \gamma = 1$
  for $A_o(I_1) = C^*(\ZZ/2\ZZ)$ and $A_u(I_1) = C^*(\ZZ)$ ; $\dim_q \gamma = 2$
  for $A_o(I_2) = C(SU_{-1}(2))$, $A_o(Q_{-1}) := C(SU(2))$ and $A_u(I_2)$.
\end{rque}

\begin{thm} \label{thm_orient_ana} Consider the quantum Cayley graph of a free
  product of universal discrete quantum groups, and assume that all directions
  have quantum dimensions different from $2$.
  \begin{enumerate}
  \item We have $p_{\ss +-}K_g = s(K_{\ss++})$.
  \item $Q_0 p_{\ss++} K_g = Q_0 K_{\ss++}$ and $(1-Q_0) p_{\ss++} K_g =
    (1-Q_0)(\id-r')(K_{\ss++})$.
  \item The subspace $\Kk_g$ is dense in $K_g$.
  \end{enumerate}
\end{thm}

\begin{dem}
  1. If $\eta = s\zeta'$ with $\zeta'\in K_{\ss++}$, put $\zeta =
  (1-r')\zeta'$. We have then $s\zeta = (1-r)\eta$, hence $\zeta\in p_{\ss++}
  K_g$ and $\eta \in p_{\ss+-}K_g$ by~\eqref{eq_JV_image}. This proves that
  $s(K_{\ss++})$ is contained in $p_{\ss+-} K_g$.

  Conversely, let $\eta\in K_{\ss+-}$, $\zeta\in K_{\ss++}$ be the projections
  of an element of $K_g$. Recall that $\eta$ is characterized by the identities
  $R_kp_k\eta = -Rs p_{\geq k}\zeta$. Since $\Ker s = Q_0 K_{\ss++}$ and $R_k$ is
  injective, we can define for each $k\in\NN$ a vector $\zeta'_k \in (1-Q_0)p_k
  K_{\ss++}$ by putting $Rs\zeta'_{k-1} = - Rs p_{\geq k}\zeta$. If the sum
  $\sum \zeta'_k$ converges to a vector $\zeta'\in K_{\ss++}$ we clearly have
  $s\zeta'=\eta$ and this shows that $p_{\ss+-}K_g \subset s(K_{\ss++})$. The
  proof of this convergence is the main issue of the Theorem and will be given
  at point~4.

  2. On the other hand the vector $\eta$ in~\eqref{eq_JV_image} lies in
  $p_{\ss+-}K_g$. Once we know that $p_{\ss+-}K_g = s(K_{\ss++})$ we can
  therefore write
  \begin{align*}
    (1-Q_0) p_{\ss++} K_g &= \{\zeta \in (1-Q_0) K_{\ss++} ~|~ 
    \exists \zeta'\in (1-Q_0) K_{\ss++}~~ (1-r)s\zeta' = s\zeta\} \\
    &= \{\zeta \in (1-Q_0) K_{\ss++} ~|~ 
    \exists \zeta'\in (1-Q_0) K_{\ss++}~~ s(1-r')\zeta' = s\zeta\} \\ 
    &= (1-Q_0)(\id-r')(K_{\ss++}),
  \end{align*}
  since $s$ is injective on $(1-Q_0) K_{\ss++}$ and vanishes on $Q_0 K_{\ss++}$.

  Finally, since $s$ vanishes on $Q_0 K_{\ss++}$ there is no obstruction
  in~(\ref{eq_JV_image}) for a vector of $Q_0 K_{\ss++}$ to belong to $p_{\ss
    ++}K_g$, hence $Q_0 p_{\ss ++} K_g = Q_0 K_{\ss++}$.

  3. If $\zeta$ is an element of $p_{\ss++}K_g$ and $\eta$ is the
  associated vector in $K_{\ss+-}$, we know \cite[proof of
  Thm.~5.3]{Vergnioux_Cayley} that $\zeta$ is the image by $p_{\ss++}$ of the
  following element $\xi\in K_g$:
  \begin{displaymath}
    \xi = \zeta - \eta - W\eta + p_{\ss--}\Theta(\eta-\zeta).
  \end{displaymath}
  Put $r'' = (1-Q_0)r'$. By injectivity of $p_{\ss++} : K_g \to K$ and the
  previous points, the equation above shows that $K_g$ is the image of the
  following continuous application:
  \begin{displaymath}
    K_{\ss++} \to K, ~~ \zeta' \mapsto (1-r'')\zeta' - (1+W)s\zeta' + 
    p_{\ss--}\Theta (s+r''-1) \zeta'.
  \end{displaymath}
  Since $\Kk_{\ss++}$ is dense in $K_{\ss++}$ and mapped into $\Kk$ by this
  application, we conclude that $\Kk_g$ is dense in $K_g$. 

  4. In the case when no direction has quantum dimension $2$, we know from
  \cite[proof of Thm.~6.5]{Vergnioux_Cayley} that $K^{-1}\|\zeta\| \leq
  \|R_k\zeta\| \leq \|\zeta\|$ for every $k$, $\zeta\in p_k K_{\ss+-}$ and a
  constant $K$ depending only on the tree. Hence we obtain from the identity
  $Rs\zeta'_{k-1} = - Rs p_{\geq k}\zeta$ the following inequality:
  \begin{equation*}
    \|s \zeta'_{k-1}\| \leq K \sum_{j=k}^\infty \|s p_j \zeta\|.
  \end{equation*}
  Now there are orthogonal subspaces $q_l P_\gamma p_j K_{\ss++}$ of $p_j
  K_{\ss++}$ on which the map $s$ is a multiple of an isometry with a norm
  known explicitly in terms of the sequence of quantum dimensions $(m_k)$
  associated to the direction $\gamma$ \cite[Lem.~6.3]{Vergnioux_Cayley}:
  \begin{displaymath}
    \|sq_lP_\gamma p_j\| = \sqrt{\frac{m_l m_{l-1}}{m_jm_{j+1}}}.
  \end{displaymath}
  In the case of $A_o$ there is only one direction $\gamma$ and $P_\gamma = \id$.
  Since $q_l P_\gamma$ commutes with all our structural maps so far
  it commutes in particular with the construction of $\zeta'_k$ from $\zeta$ and
  we obtain:
  \begin{equation*}
    \|q_lP_\gamma\zeta'_{k-1}\| \leq K \sum_{j=k}^\infty 
    \sqrt{\frac{m_{k-1}m_k}{m_jm_{j+1}}} \|p_jq_lP_\gamma\zeta\|. 
  \end{equation*}
  In the case when the direction $\gamma$ has quantum dimension different from
  $2$, its sequence of dimensions has the form $m_{k-1} = (a^k - a^{-k}) / (a -
  a^{-1})$ for some $a>1$ and it is easy to check that one has then
  $\frac{m_{k-1}}{m_j} \leq a^{-(j-k)}$ for all $j\geq k-1$. This allows to write:
  \begin{tabbing}
    \hspace{3em} \= \hspace{1.3em} $\ds\|q_lP_\gamma\zeta'_{k-1}\|^2$ \= $\ds \leq K^2
    \Big(\sum_{j=k}^\infty a^{k-j} \|p_jq_lP_\gamma\zeta\|\Big)^2$ \\
    \> \> $\ds \leq K^2 \sum_{j=k}^\infty a^{k-j} 
    \sum_{j=k}^\infty a^{k-j} \|p_jq_lP_\gamma\zeta\|^2 $ \\
    \> \> $\ds = \frac{K^2}{1-a^{-1}}~ \sum_{j=k}^\infty a^{k-j}
    \|p_jq_lP_\gamma\zeta\|^2$,~~~ hence \\
    \> $\ds \sum_{k=1}^\infty \|q_lP_\gamma\zeta'_{k-1}\|^2$ \> 
    $\ds \leq \frac{K^2}{1-a^{-1}}~
    \sum_{1\leq k\leq j} a^{k-j} \|p_jq_lP_\gamma\zeta\|^2$ \\
    \> \> $\ds \leq \frac{K^2}{(1-a^{-1})^2}~ \sum_{j=1}^\infty 
    \|p_jq_lP_\gamma\zeta\|^2 
    = \frac{K^2 \|q_lP_\gamma\zeta\|^2} {(1-a^{-1})^2}.$
  \end{tabbing}
  Finally we sum over $l$, take the smallest $a$ when $\gamma$ varies, sum over
  $\gamma$ and apply the ``cut-and-paste'' process
  \cite[Rem.~6.4(2)]{Vergnioux_Cayley} to recover the whole of the tree,
  i.e. the whole of $(\zeta'_k)$ and $\zeta$. Since the projections $P_\gamma$
  on the one hand, and $q_l$ on the other hand, are mutually orthogonal, the
  inequality above shows that $(\sum \zeta'_k)$ converges.
\end{dem}

\section{Path cocycles in quantum Cayley trees} 
\label{sec_cayley}

In this section we go on with the search for paths in quantum Cayley graphs of
universal quantum groups, by proving injectivity results for $\B : K_g \to
H$. We will then examine the triviality of the cocycles obtained in this way,
distinguishing the ``orthogonal case'', and the remaining ``general'' cases of
free products of unitary and orthogonal universal discrete quantum groups.

\bigskip

Recall that a classical graph is a tree \iff $\B : \Kk_g \to \Ker \epsilon
\subset \Hh$ is bijective. Hence our study has a geometrical interpretation,
namely the question whether the quantum Cayley graphs of universal quantum
groups are trees. We have already seen in \cite{Vergnioux_Cayley} that the
ascending orientation of our quantum Cayley graphs behaves like in rooted trees,
e.g. $\B : K_{\ss++} \to H$ is injective.

One cannot expect the map $\B : K_g \to H$ to be injective at the $\ell^2$
level in general: it is already not the case for non-abelian free groups because
of exponential growth. Surprisingly, the results of the previous section yield a
strong ``$\ell^2$-injectivity'' result for the ``purely quantum'' part
$(1-Q_0)K_g$ of the space of antisymmetrical edges.

\begin{prp} \label{prp_q_path} Consider the quantum Cayley graph of a free
  product of universal discrete quantum groups. If all directions have quantum
  dimensions different from $2$, then the restriction of $\B$ to $Q_0 \Kk_g
  \oplus (1-Q_0) K_g$ is injective. Moreover $\B (1-Q_0) K_g$ is a closed
  subspace of $H$.
\end{prp}

\begin{dem}
  Since $\B(1-Q_0)K_g \subset (1-Q_0)H$ by Lemma~\ref{lem_quasiclass}, it
  suffices to show separately that $\B$ is injective on $(1-Q_0)K_g$, and that
  $Q_0\B$ is injective on $Q_0 \Kk_g$. We start with the first assertion and
  prove at the same time that $\B (1-Q_0) K_g$ is closed.

  Observe that $(1-Q_0)K_g \subset K_g$ because $[Q_0, \Theta] = 0$.  Since the
  restrictions of $\B$ and $(\B-\O)$ coincide on $K_g$ up to a constant, and by
  taking the adjoint, it is sufficient to show that the inclusion
  \begin{displaymath}
    (1-Q_0)(\B-\O)^*(H) \subset (1-Q_0)K_g
  \end{displaymath}
  given by Lemma~\ref{lem_B-O} is in fact an equality.  It is moreover
  equivalent to show that $(1-Q_0)p_{\ss++}(\B-\O)^*(H) = (1-Q_0)p_{\ss++} K_g$
  because $p_{\ss++} : K_g \to K_{\ss++}$ is injective with closed range.

  Therefore we consider the following operator:
  \begin{align*}
    p_{\ss++} (\B-\O)^* &= p_{\ss++}(1-\Theta)\B^* = 
    p_{\ss++}\B^* - p_{\ss++}\Theta p_{\ss--}\B^* \\
    &= (1 - p_{\ss++}\Theta p_{\ss--}\Phi^*)p_{\ss++}\B^* + p_{\ss++}\Theta
    p_{\ss--} \B^* p_0.
  \end{align*}
  According to this expression, and because $\B(K_{\ss++}) = (1-p_0)H$ by
  \cite[Prop.~4.7]{Vergnioux_Cayley}, we have $p_{\ss++}(\B-\O)^*(H) \supset (1
  - r') (K_{\ss++})$, and we obtain by Theorem~\ref{thm_orient_ana} the 
  inclusion we are looking for:
  \begin{displaymath}
    (1-Q_0)p_{\ss++}(\B-\O)^*(H) \supset (1-Q_0)p_{\ss++} (K_g).
  \end{displaymath}

  Now let us prove that $Q_0 \B$ is injective on $Q_0 \Kk_g$, by
  contradiction. Fix $\zeta \in Q_0 \Kk_g$ such that $\zeta \neq 0$ and $Q_0 \B
  \zeta = 0$. For each edge $w$ of the classical Cayley tree, $w =
  (\alpha,\beta)$ with direction $\gamma$, we define a projection $p_w \in B(K)$
  by putting $p_w = (p_\alpha\tens p_\gamma) p_{\ss++}$ if $w$ is ascending,
  $p_w = (p_\alpha\tens p_\gamma) p_{\ss--}$ otherwise. In our situation it is
  known that $\{p_w\}$ is a partition of unity in $B(K)$
  \cite[Prop.~4.5]{Vergnioux_Cayley}, which clearly stabilizes $\Kk$. 

  In particular the support $\Supp \zeta = \{w ~|~ p_w \zeta \neq 0\}$ of
  $\zeta$ consists of edges of a disjoint union of finite subtrees of the
  classical Cayley tree. Note also that $\Supp \zeta$ is stable under
  the classical reversing map, since $\Theta \zeta = -\zeta$ and $\zeta \in
  K_{\ss++} \oplus K_{\ss--}$: for $w = (\alpha,\beta)$ in the support of
  $\zeta$, e.g. ascending with direction $\gamma$, we have by
  \cite[Prop.~3.7]{Vergnioux_Cayley}:
  \begin{align*}
    0\neq \Theta (p_\alpha\tens p_\gamma) p_{\ss+\jok} \zeta &=
    \hat\delta(p_\alpha)(\id\tens p_{\bar\gamma}) p_{\ss\jok-}\Theta \zeta
    = \hat\delta(p_\alpha)(\id\tens p_{\bar\gamma}) p_{\ss\jok-} \zeta \\
    &= \hat\delta(p_\alpha)(\id\tens p_{\bar\gamma}) p_{\ss--} \zeta =
    (p_\beta\tens p_{\bar\gamma}) p_{\ss--} \zeta,
  \end{align*}
  hence $\bar w = (\beta,\alpha)$ belongs to $\Supp \zeta$.

  We fix one of the subtrees in $\Supp \zeta$, and we choose one of its terminal
  edges $w_1$ wich is not the closest one to the origin. In particular, $w_1$
  must be ascending and the target of any other $w \in \Supp \zeta$ is different
  from the target $\beta_1$ of $w_1$. Now $\B p_{w_1} \zeta$ is a non-zero
  element of $p_{\beta_1}H$, and it equals $Q_0 \B p_{w_1} \zeta$ by the last
  point of Lemma~\ref{lem_quasiclass}. For any other edge $w$ with target
  $\beta$ the vector $Q_0 \B p_w \zeta$ lies in $Q_0 p_{\beta} H$, which is
  orthogonal to $p_{\beta_1}H$. Hence we obtain $Q_0 \B \zeta \neq 0$, a
  contradiction.
\end{dem}

\begin{crl} \label{crl_path_tree} Consider a quantum Cayley graph as above where
  all directions have quantum dimension different from $2$. Then there exists a
  unique path cocycle $c_g : \Ss \to \Kk'_g$.
\end{crl}

\begin{dem}
  Since $\Theta^2 = \id$ on $Q_0 K$, we have $p_g = (\id-\Theta)/2$ on $Q_0 K$,
  hence $Q_0 \Kk'_g = Q_0 \Kk_g$. As a result, we see that $\Kk'_g \subset Q_0
  \Kk_g \oplus (1-Q_0)K_g$, and we can apply Corollary~\ref{crl_path_prop} with
  $K'_g = \Kk'_g$.
\end{dem}

\subsection{The orthogonal case}
\label{sec_ortho_path}

In this section we recall the definition of the classical projections $q_0$ and
we describe the ``classical subgraph'' of the quantum Cayley graph obtained in
this way. This description holds generally for free products of orthogonal and
unitary universal discrete quantum groups, and it is strongly related to the
structure of the underlying classical Cayley graph.

Only then we will restrict to the case where there is a unique factor $A_o(Q)$
in the free product quantum group under consideration. In this case the
projections $q_0$, $Q_0$ coincide and, using Proposition~\ref{prp_q_path}, we
will be able to deduce that the target operator is in fact invertible from $K_g$
to $H$. In particular, the path cocycle obtained previously is trivial.

\begin{rque}{Reminder} \label{rem_ql} The classical subspaces $q_0H = H_0$,
  $q_0K = K_0$ can be introduced as subspaces of fixed vectors for the
  left-right representations of $\hat S$ on $H$, $K$. More precisely, we have
  two (resp. four) commuting representations of $\hat S$ on $H$ (resp. $K$)
  given by $\hat\pi_2(x\tens x') = x(Ux'U) \in B(H)$ and $\hat\pi_4(x\tens
  y\tens y'\tens x') = (x\tens y)(Ux'U\tens Uy'U)$, and the left-right
  representations above are $\hat\pi_2\rond\hat\delta$ and
  $\hat\pi_4\rond\hat\delta^3$. Hence we have $q_0 = \hat\pi_2\hat\delta(p_0)$
  on $H$ and $q_0 = \hat\pi_4\hat\delta^3(p_0)$ on $K$ --- more generally one
  defines in fact a new partition of unity $\{q_l\}$ by replacing $p_0$ with
  $p_{2l}$ in these formulae.

  With this definition it is easy to check that the projections $q_l$ commute to
  the structural maps $\B$, $\Theta$ \cite[Prop.~3.7]{Vergnioux_Cayley} and,
  directly from the definitions, to the partitions of unity $\{p_n\}$,
  $\{p_{\ss\jok+}, p_{\ss\jok-}\}$, $\{p_{\ss+\jok}, p_{\ss-\jok}\}$. In
  particular, $q_0 \in B(H)$ and $q_0 \in B(K)$ can be viewed as projections
  onto a ``classical subtree'' $(H_0, K_0)$. Notice however that the subspaces
  $H_0$, $K_0$ are not stable under the action of $\Ss$ in the quantum case. In
  the classical case $\hat S$ is commutative and we have thus $H_0 = H$, $K_0 =
  K$.

  Moreover we have in the general case $p_{\ss+-} K_0 = p_{\ss -+} K_0 = 0$
  hence $\Theta$ restricts to an involution on $K_0$
  \cite[Prop.~5.1]{Vergnioux_Cayley}, and $q_0K \subset Q_0K$. In the case of
  $A_o(I_n)$, this inclusion is in fact an equality. As a matter of fact by
  \cite[Lemma~6.3]{Vergnioux_Cayley} we have $q_0 K_{\ss ++} = \Ker s = Q_0
  K_{\ss++}$. By applying $\Theta$ we also get $q_0 K_{\ss --} = Q_0 K_{\ss--}$,
  hence finally $q_0 = Q_0$.
\end{rque}

Let us now describe this classical subtree $(H_0, K_0)$. We identify $H$ not
only with the GNS space of $h$, but also with the GNS space of the left Haar
weight $\hat h_L$ of $\hat S$. For each vertex $\alpha \in \Irr\Cc$ of the
classical Cayley graph, the subspace $p_\alpha H_0 \subset p_\alpha H$
corresponds to the inclusion $1_\Cc \subset \bar\alpha \tens \alpha$ of
representations of $\hat S$, hence it is $1$-dimensional and spanned by the
vector $\xi_\alpha \in p_\alpha H$ which is the GNS image of $p_\alpha \in \hat
S$. Note that the vector $\xi_\alpha$ associated to the trivial representation
is $\xi_0 = \Lambda_h(1)$.

Similarly, we define $\xi_{(\alpha,\beta)}$ for any edge $(\alpha,\beta)$ of the
classical Cayley graph as the GNS image in $K = H\tens p_1H$ of
$\hat\delta(p_\beta)(p_\alpha\tens p_1)$, which spans the subspace
$\hat\delta(p_\beta)(p_\alpha\tens p_1) K_0$ and is mapped into $p_\alpha H_0$
(resp. $p_\beta H_0$) by $\O$ (resp. $\B$). Denote by $\gamma$ the direction of
the edge from $\alpha$ to $\beta$ in the classical Cayley graph, so that
$\hat\delta(p_\beta)(p_\alpha\tens p_1) = \hat\delta(p_\beta)(p_\alpha\tens
p_\gamma)$. The norms of the vectors $\xi_\alpha$ and $\xi_{(\alpha,\beta)}$ are
easy to compute from the known formulae for $\hat h_L$, see
e.g. \cite[(2.13)]{PodlesWoro_Lorentz}: 
\begin{align*}
  \|\xi_\alpha\|^2 &= m_\alpha \Tr_\alpha (Fp_\alpha) = m_\alpha^2 \text{~~~and}
\end{align*}
\begin{align*}
  \|\xi_{(\alpha,\beta)}\|^2 &= m_\alpha m_\gamma (\Tr_\alpha\tens \Tr_\gamma)
  ((F\tens F)\hat\delta(p_\beta)) \\ &= m_\alpha m_\gamma \Tr_\beta (Fp_\beta)
  = m_\alpha m_\beta m_\gamma.
\end{align*}
Here $m_\alpha = \Tr_\alpha(Fp_\alpha)$ is the quantum dimension of $\alpha$.

Since $q_0 K_{\ss+-} = q_0 K_{\ss-+} = \{0\}$, the operator $\Theta$ maps
$\hat\delta(p_\beta)(p_\alpha\tens p_1) K_0$ isometrically to
$\hat\delta(p_\alpha) (p_\beta\tens p_1) K_0$, hence
$\Theta(\xi_{(\alpha,\beta)}) = \xi_{(\beta,\alpha)}$ --- recall that the
direction of $(\beta,\alpha)$ is $\bar\gamma$, and $m_{\bar\gamma} = m_\gamma$.
In particular, the vectors
\begin{equation*}
  \tilde \xi_{\alpha\wedge\beta} = \frac{\xi_{(\alpha,\beta)} -
    \xi_{(\beta,\alpha)}} {\sqrt{2\,m_\alpha m_\beta m_\gamma}} \text{~~~and~~~}
  \tilde \xi_\alpha = \frac 1{m_\alpha}~ \xi_\alpha,
\end{equation*}
where $(\alpha,\beta)$ runs over the classical ascending edges, and $\alpha$,
over the classical vertices, form hilbertian bases of $q_0 K_g$, $q_0 H$
respectively. 

On the other hand, the norm of $\B$ on the subspaces
$\hat\delta(p_\beta)(p_\alpha\tens p_\gamma) K$ is known from
\cite[Rem.~4.8]{Vergnioux_Cayley} and we deduce
\begin{equation*}
  \B (\xi_{(\alpha,\beta)}) = \frac{m_\alpha m_\gamma}{m_\beta}~ \xi_\beta, ~~~
  \B (\tilde \xi_{\alpha\wedge\beta}) = \sqrt{\frac{m_\gamma}2}~ \left(
  \sqrt{\frac{m_\alpha}{m_\beta}}~ \tilde\xi_\beta - 
  \sqrt{\frac{m_\beta}{m_\alpha}}~ \tilde\xi_\alpha\right).
\end{equation*}
We see that $(H_0, K_0, \Theta, \B)$ is a hilbertian version of the classical
Cayley tree where the relation ``$\alpha$ is an endpoint of $w$'' comes with
weights depending on the quantum dimensions.

As a result we can find unique paths in $q_0 \Kk_g$ in the following sense: for
each vertex $\alpha$ of the classical Cayley tree, the vector $\tilde\xi_\alpha/
m_\alpha - \xi_0$ has a unique preimage by $\B$ in $\Kk_g$. Denoting by
$(1=\alpha_0, \alpha_1, \ldots, \alpha_{n-1}, \alpha_n=\alpha)$ the geodesic
from $1_\Cc$ to $\alpha$ in the classical Cayley tree, and by $\gamma_i$ the
direction from $\alpha_i$ to $\alpha_{i+1}$, this preimage is given by
\begin{displaymath}
  \zeta_\alpha =  \sum_{i=0}^{n-1} \ts\sqrt{\frac 2{m_{\gamma_i}}}\ds ~
  \frac {\tilde\xi_{\alpha_i\wedge\alpha_{i+1}}}
  {\sqrt{m_{\alpha_i} m_{\alpha_{i+1}}}}.
\end{displaymath}
Moreover, for each infinite geodesic $\infty = (1, \alpha_1, \ldots, \alpha_n,
\ldots)$ in the classical Cayley tree such that $(\sum m_{\alpha_i}^{-2})$
converges we get a preimage of $\xi_0$ :
\begin{equation}\label{eq_fixed}
  \zeta_{\infty} = 
  \sum_{i=0}^\infty \ts\sqrt{\frac 2{m_{\gamma_i}}}\ds~ 
  \frac {\tilde\xi_{\alpha_i\wedge\alpha_{i+1}}}
  {\sqrt{m_{\alpha_i} m_{\alpha_{i+1}}}}.
\end{equation}
In the case of $A_o(Q)$ there is a unique geodesic to $\infty$ and the previous
vector is the fixed vector $\xi_g$ for the path cocyle $c_g$ of the next
Theorem:

\begin{thm} \label{thm_path_Ao} Consider the quantum Cayley graph of $A_o(Q)$,
  with $Q^*Q \notin \RR I_2$. Then the operator $\B : K_g \to H$ is invertible.
  As a result there exists a unique path cocycle $c_g : \Ss \to \Kk'_g$, and it
  is trivial.
\end{thm}

\begin{dem}
  The classical Cayley graph of $A_o(Q)$ is the half-line with vertices
  $\alpha_k$ at integers $k\in\NN$. We put $\tilde\xi_k = \tilde\xi_{\alpha_k}$
  and $\tilde\xi_{k\wedge k+1} = \tilde\xi_{\alpha_k\wedge\alpha_{k+1}}$. If
  $\B\zeta=0$ with $\zeta = \sum\lambda_k\tilde\xi_k$, the equality $p_0\B\zeta
  = 0$ reads $\lambda_0 = 0$ and an easy induction yields $\lambda_k = 0$ for
  all $k$: $\B$ is injective on $q_0 K_g$. Define conversely a linear map
  $\B^{-1} : q_0\Hh \to q_0K_g$ by
  \begin{displaymath}
    \B^{-1}(\tilde\xi_k) = \ts-m_k\sqrt{\frac 2{m_1}}\ds~ \sum_{i=k}^\infty
    \frac {\tilde \xi_{i\wedge i+1}}{\sqrt{m_i m_{i+1}}}.
  \end{displaymath}
  Clearly $\B \B^{-1} = \id$, $\B^{-1} {\B}_{|q_0\Kk_g} = \id$ and it remains to
  check that $\B^{-1}$ extends to a bounded operator on $q_0 H$. Using the
  asymptotics $m_k \sim C a^k$ with $a > 1$ we get
  \begin{displaymath}
    0 \leq \Big(\B^{-1}(\tilde\xi_k) \Big| \B^{-1}(\tilde\xi_l) \Big) =
    \frac 2{m_1} \sum_{i = \max(k,l)}^\infty \frac {m_km_l}{m_im_{i+1}}
    \leq D a^{-|k-l|}
  \end{displaymath}
  for some constant $D$, and the matrix $(a^{-|k-l|})_{k,l}$ is easily seen to
  be bounded.

  This proves that $\B : q_0 K_g \to q_0 H$ is invertible.  Since $q_0 = Q_0$,
  Proposition~\ref{prp_q_path} shows that $\B : (1-q_0)K_g \to (1-q_0)H$ is
  injective with closed range. Moreover, by Lemma~\ref{lem_B-O} the subspace
  $\B(K_g)$ contains $\Ker \epsilon$, which is dense in $H$. Hence $\B : K_g \to
  H$ is invertible. In particular Corollary~\ref{crl_path_prop} applies with
  $K'_g = K_g$. Finally $\xi_0 = \Lambda_h(1)$ is tautologically a fixed vector
  for $c_0$, and the unique vector $\xi_g \in K_g$ such that $\B\xi_g = \xi_0$
  is a fixed vector for $c_g$.
\end{dem}

\subsection{The general case}

We come back to the general case of free products of orthogonal and unitary
universal quantum groups. We will prove that the path cocycle obtained at
Corollary~\ref{crl_path_tree} is not trivial --- except in the cases considered
in the previous section. The strategy is as follows: we will exhibit a
particularly simple infinite geodesic in the quasi-classical subtree, compute
the values of the path cocycle along this geodesic, and realize that they are
unbounded.

Notice that the computations of the previous section easily yield the existence
in general of ``weak'' path cocycles which are trivial, i.e. of trivial cocycles
$c : \Ss \to K_g$ such that $\B \rond c_g = c_0$: for each infinite geodesic
starting from $1$ in the classical Cayley graph, we get a vector $\zeta_\infty
\in K_g$ such that $\B \zeta_\infty = \xi_0$, and the trivial cocycle with fixed
vector $\zeta_\infty$ is clearly a ``weak'' path cocycle. However such a path
cocycle does not take its values in $\Kk_g$ nor $\Kk'_g$ in the non-orthogonal
case. Notice also in the non-orthogonal case that we can choose two different
geodesics $\infty$, $\infty'$, and obtain in this way a non-zero vector $\zeta =
\zeta_\infty - \zeta_{\infty'} \in K_g$ such that $\B \zeta = 0$.

\bigskip

Of course we have to exclude now the case when the quantum group under
consideration consists of a single copy of $A_o(Q)$, already studied in
Section~\ref{sec_ortho_path}. For the sake of simplicity we assume that our
quantum group contains a copy of some $A_u(Q)$. We denote by $\gamma \in
\Irr\Cc$ the corresponding fundamental corepresentation and we recall that the
tensor powers $\gamma^{\tens n}$ are irreducible and denoted $\gamma^n \in
\Irr\Cc$. In the case where the quantum group does not contain any $A_u(Q)$, but
instead two free copies of $A_o(Q_1)$, $A_o(Q_2)$ with fundamental
corepresentations $\gamma_1$, $\gamma_2$, the corepresentations $\gamma^n$
should be replaced by alternating irreducible tensor products $\gamma_1\gamma_2
\gamma_1 \cdots \gamma_k$.

Put $H_1 = p_{\gamma} H$, denote $\tilde\xi_1 = \tilde\xi_\gamma \in H_1$ the
normalized GNS image of $p_\gamma \in \hat S$ already introduced in
Section~\ref{sec_ortho_path}, and define $H_1^\circ = \tilde\xi_1^\bot \cap
H_1$. One can check that $\epsilon(\tilde\xi_1) = m_1 := \dim_q\gamma$. Due to
the irreducibility mentioned above, we have natural Hilbert space
identifications $p_{\gamma^n} H = H_1^{\tens n}$, $(p_{\gamma^n} \tens p_\gamma)
K = (p_{\gamma^n} \tens p_\gamma) K_{\ss++} = H_1^{\tens n} \tens H_1$.  For any
$0\leq l\leq n$ we define the ``length $l$'' subspaces:
\begin{gather*}
  q_l H_1^{\tens n} = H_1^{\tens l-1}\tens H_1^\circ\tens \tilde\xi_1^{\tens n-l}, 
\end{gather*}
with the convention that $q_0 H_1^{\tens n} = \CC \tilde\xi_1^{\tens n}$. This
defines projections $q_l \in B(p_{\gamma^n}H)$ and $q_l \in
B((p_{\gamma^{n-1}}\tens p_\gamma)K)$ for any $n\in\NN^*$ and $0\leq l\leq n$.

This is of course compatible with the notation $q_0$ used in
Section~\ref{sec_ortho_path}, and more generally this gives an explicit
description of the projections $q_l = \hat\pi_2\hat\delta(p_{2l})$, $q_l =
\hat\pi_4 \hat\delta^3(p_{2l})$ of \cite[p.~125]{Vergnioux_Cayley} on the
(particularly simple) subspaces $p_{\gamma^n}H$, $(p_{\gamma^n}\tens
p_\gamma)K$. In particular it is known that $\B$, $\O$, $\Theta$ intertwine the
relevant projections $q_l$. Concretely, $\B : (p_{\gamma^n}\tens p_\gamma)K \to
p_{\gamma^{n+1}}H$ is just the identification $H_1^n\tens H_1 = H_1^{\tens
  n+1}$, the operator $\O = \id\tens\epsilon$ maps $\zeta\tens\tilde \xi_1 \in
q_l(p_{\gamma^n}\tens p_\gamma)K$ to $m_1 \zeta$ if $l\leq n$, and
$q_{n+1}(p_{\gamma^n}\tens p_\gamma) K$ to $0$.

\bigskip

Now we fix $n$, $l$, and $\zeta \in q_l p_{\gamma^n} H$. As in the proof of
Lemma \ref{lem_B-O}, we can easily construct an element $\eta \in \Kk_{\ss++}$
such that $(\B-\O)\eta = \zeta - \epsilon(\zeta) \xi_0$: we put $\eta = \eta_1 +
\eta_2 + \cdots + \eta_n$ with $\eta_i \in (p_{\gamma^{n-i}}\tens
p_\gamma)K_{\ss++}$ uniquely determined by the identity $\B\eta_i =
\O\eta_{i-1}$, $\B\eta_1 = \zeta$. By the remarks in the preceding paragraph we
see that $\eta_i$ vanishes for $i > n-l+1$, and $\|\eta_i\| =
m_1^{i-1}\|\zeta\|$ for $i\leq n-l+1$.

To obtain a ``path from the origin to $\zeta$'' in $\Kk'_g$, it remains to
project $\eta$ onto $K_g$. This is easy because $\eta$ lies ``almost entirely''
in $Q_0 K_{\ss++}$. More precisely, we assert that $\eta_i$ belongs to $Q_0
K_{\ss++}$ for $i \leq n-l$, and that $\eta_{n-l+1}$ belongs to
$(1-Q_0)K_{\ss++}$.

As a matter of fact, the subspace $(p_{\gamma^{n+1}}\tens p_{\bar\gamma})
K_{\ss+-}$ is isomorphic as a left-right representation of $\hat S$ to
$\gamma^{n+1}\bar\gamma \tens \bar\gamma^n =
\gamma^{n+1}\bar\gamma^{n+1}$, and in particular it is irreducible of length
$n+1$, so that $q_l (p_{\gamma^{n+1}}\tens p_{\bar\gamma}) K_{\ss+-} = 0$ if
$l\leq n$ \cite[p.~125]{Vergnioux_Cayley}. Since we have
$\Theta(p_{\gamma^n}\tens p_\gamma) K_{\ss++} \subset (p_{\gamma^{n+1}}\tens
p_{\bar\gamma}) K_{\ss\jok-}$ , this implies that $q_l(p_{\gamma^n}\tens
p_\gamma) K_{\ss++} \subset Q_0 K_{\ss++}$ for $l\leq n$. Similarly,
$q_{n+1}(p_{\gamma^{n+1}}\tens p_{\bar\gamma}) K_{\ss--}$ vanishes, so that
$q_{n+1} (p_{\gamma^n}\tens p_\gamma)K_{\ss++} \subset (1-Q_0) K_{\ss++}$.

From these remarks we can conclude that the orthogonal projection of $\eta_i$
onto $K_g$ simply equals $(\eta_i-\Theta(\eta_i))/2$ and is orthogonal to the one of
$\eta_{n-l+1}$, for each $i < n-l+1$. Summing over these $i$'s we obtain
\begin{equation} \label{eq_unit_path_min}
  \|c_g(\zeta)\|^2 = \|p_g(\eta)\|^2 \geq \ts \frac 12 
  \sum_{i=1}^{n-l} m_1^{2(i-1)} \|\zeta\|^2
  \geq C m_1^{2(n-l)} \|\zeta\|^2
\end{equation}
for some constant $C>0$ depending only on $m_1$. This should be put in contrast
with the situation in $A_o(Q)$, where we have $\|c_g(\zeta)\| \leq C \|\zeta\|$
for all $\zeta$ such that $\epsilon(\zeta) = 0$, by
Theorem~\ref{thm_path_Ao}. Note however that this does not directly mean that
the path cocycle is ``unbounded''. The precise statement is indeed the following
one:

\begin{prp} \label{prp_unitary_nontriv}
  Let $\gamma$ be the fundamental corepresentation of $A_u(I_N)$ and consider
  the corepresentations $\gamma^n \in B(H_{\gamma^n})\tens S$. Then the norm of
  $C_n := (\id\tens c_g)(\gamma^n)$ as an element of $B(H_{\gamma^n},
  H_{\gamma^n}\tens K_g)$ is bounded below by $C \sqrt{n+1}$ for some constant
  $C>0$ depending only on $m_1$. In particular the path cocycle $c_g : \Ss \to
  \Kk'_g$ is not trivial in the non-orthogonal case.
\end{prp}

\begin{dem}
  We fix an orthonormal basis $(e_i)$ of the space $H_\gamma$ of the
  corepresentation $\gamma$, and we denote by $u_{ij}$ the corresponding
  generators of $A_u(Q)$. If $\underline i = (i_1, \ldots, i_n)$, $\underline k
  = (j_n, \ldots, j_n)$ are multi-indices, we put $e_{\underline i} =
  e_{i_1}\tens \cdots \tens e_{i_n} \in H_{\gamma^n}$ and $u_{\underline i
    \underline k} = u_{i_1k_1} \cdots u_{i_nk_n} \in \Ss$. Via the GNS map $\Ss
  \to H$ the vector $u_{\underline i \underline k}$ identifies with an element
  of $p_{\gamma^n}H$.

  We have by definition and Inequality~\eqref{eq_unit_path_min}:
  \begin{align*}
    \|C_n(e_{\underline i})\|^2 &= \|\ts\sum_{\underline k} 
    e_{\underline k} \tens c_g(u_{\underline i \underline k})\|^2 
    = \ts\sum_{\underline k} \|c_g(u_{\underline i \underline k})\|^2 \\
    & \geq C \ts\sum_{\underline k}\sum_{l=0}^n m_1^{2(n-l)} 
    \|q_l(u_{\underline i \underline k})\|^2.
  \end{align*}
  Observe indeed that $c_g$ maps the respective subspaces $q_l \Hh$ to the
  subspaces $q_l K_g$, which are mutually orthogonal. 

  It remains to find a lower bound for $\|q_l(u_{\underline i \underline
    k})\|$. This is where we use the additional ``unimodularity'' assumption
  (i.e. $Q = I_N$) to simplify computations. In this case we can indeed identify
  $p_{\gamma^n} H$ with $B(H_{\gamma^n}) \simeq B(H_\gamma)^{\tens n}$ equipped
  with the normalized Hilbert-Schmidt norm $\|a\|^2_{HS} = \Tr (a^*a)/m_1^n$,
  by sending $u_{\underline i\underline k}$ to $e_{\underline i} e_{\underline
    k}^*$. By definition, $q_l p_{\gamma^n} H$ corresponds then to the subspace
  of applications of the form $a_1 \tens \cdots a_l \tens \id\tens \cdots \tens
  \id$ with $a_i \in B(H_\gamma)$ and $\Tr(a_l) = 0$.

  In particular, if $k_l \neq i_l$ and $k_p = i_p$ for $p>l$ we have in this
  identification
  \begin{displaymath}
    q_l(e_{\underline i} e_{\underline k}^*) = m_1^{-(n-l)}~ (e_{i_1}e_{k_1}^*) \tens \cdots 
    \tens (e_{i_l}e_{k_l}^*) \tens \id \tens \cdots\tens \id.
  \end{displaymath}
  Hence we obtain $\|q_l(u_{\underline i \underline k})\|^2 = m_1^{-2(n-l)}
  m_1^{n-l} / m_1^n = m_1^{-2n+l}$. Observing moreover that there are at least
  $m_1^{l-1}$ multi-indices $\underline k$ satisfying the conditions above for a
  fixed $\underline i$, our first inequality yields
  \begin{displaymath}
    \|C_n(e_{\underline i})\|^2 \geq C \ts\sum_{l=0}^n m_1^{2(n-l)} m_1^{l-1} m_1^{-2n+l}
    = (n+1)\frac {C}{m_1}.
  \end{displaymath}

  Finally, if $c_g$ was trivial with fixed vector $\xi_g \in K_g$, we would have
  $C_n(\zeta) = \gamma^n(\zeta\tens\xi_g)$, hence $\|C_n(\zeta)\|
  \leq \|\zeta\| \times \|\xi_g\|$, for all $n$ and $\zeta \in H_{\gamma^n}$.
\end{dem}

To conclude this section, let us notice that the path cocycle $c_g$ considered
in the previous Proposition is unbounded, but at the same time it can be shown
not to be proper. To see this, one considers the subspaces of $H$, $K$
corresponding to the infinite geodesic $(1, \gamma, \gamma\bar\gamma,
\gamma\bar\gamma\gamma, \ldots)$ in the classical Cayley graph of $\Aa_u(Q)$ ---
they are given by the projection $P_\gamma$ already used in the proof of
Theorem~\ref{thm_orient_ana} and introduced in
\cite[p.~125]{Vergnioux_Cayley}. Then the results of the orthogonal case apply
to these subspaces, basically because $q_0 P_\gamma= Q_0 P_\gamma$. In
particular one can show that the restriction of $c_g$ to this particular
geodesic is bounded.

\section{Application to the first $L^2$-cohomology}
\label{sec_cohomology}

We begin this section with some notation. Let $\Ss$ be the dense Hopf algebra
associated with a discrete quantum group.  For any representation $\pi : \Ss \to
L(X)$ of $\Ss$ on a complex vector space $X$, we endow $X$ with the trivial
right $\Ss$-module structure given by $\epsilon$, and we denote by $H^1(\Ss,X)$
the first Hochschild cohomology group of $\Ss$ with coefficients in $X$. Recall
that we have $H^1(\Ss, X) = \Der (\Ss, X) / \Inn (\Ss, X)$, where $\Der(\Ss, X)$
coincides with the group of $\pi$-cocycles from $\Ss$ to $X$, and \linebreak $\Inn(\Ss,X)$
is the subgroup of trivial cocycles. 

The coefficient modules of interest for this section will be $H$, via the
regular representation $\lambda : \Ss \to B(H)$, and the von Neumann algebra $M
:= \lambda(\Ss)'' \subset B(H)$, via left multiplication by elements of
$\Ss$. 

In the previous sections we have studied one particular cocycle $c$ on $\Ss$,
namely the path cocycle $c_g$ with values in $K$. In the case of the universal
quantum groups $A_o(Q)$, we have seen that $c_g$ is trivial as an element of
$H^1(\Ss, K)$. In this section we use the vanishing of $c_g$ to prove that many
more cocycles are trivial: in fact the whole $L^2$-cohomology group $H^1(\Ss,
M)$ vanishes in the case of $A_o(I_n)$. The heuristic reason is that $c_g$,
being a push-back of the trivial cocycle through the restricted multiplication
map $\B$, is sufficiently universal amongst $L^2$-cocycles.

In fact a slight strengthening of the vanishing of $c_g$ will be needed: we will
need to know that the fixed vector $\xi_g \in K = H\tens p_1H$ lies in
$\Lambda_h(M)\tens p_1H$.

\begin{thm} \label{thm_vanishing} Let $\Ss$, $M$, $H$ be the Hopf algebra, the
  von Neumann algebra and the Hilbert space associated with a unimodular
  discrete quantum group. Using the terminology of
  Definition~\ref{df_path_cocycle} we assume that there exists a path cocycle
  $c_g : \Ss \to K$ which is trivial, and which admits moreover a fixed vector
  $\xi_g$ lying in $\Lambda_h(M)\tens p_1H \subset K$. Then we have
  $H^1(\Ss,M)=0$.
\end{thm}

\begin{dem}
  We proceed in 4 steps.

  1. We observe that for any linear map $c : \Ss \to M$ the following formula,
  making use of the right $M$-module structure of $H$, defines a bounded map
  $m_c : K \to H$:
  \begin{displaymath}
    m_c(\zeta\tens\Lambda(x)) = \zeta \cdot c(x).
  \end{displaymath}
  The main reason is that in this formula only the values of $c$ on the
  finite-dimensional subspace $\Ss_1 \subset \Ss$ corresponding to $p_1 H\subset
  \Hh$ play a role.

  More precisely, recall that the Haar state $h$ is a trace in the unimodular
  case, and denote by $\rho(x)\in B(H)$ the action by right multiplication of
  $x\in M$. Fix an orthonormal basis $(\Lambda(x_i))_{1\leq i\leq N}$ of the
  space $p_1H$. Any $\eta \in K$ can be decomposed into a sum $\sum
  \zeta_i\tens\Lambda(x_i)$ and we have then
  \begin{displaymath}
    \|m_c(\eta)\| = \|\ts\sum \rho(c(x_i))\zeta_i\| 
    \leq \sum \|c(x_i)\|_M \|\zeta_i\|_H \leq C \sqrt N \|\eta\|,
  \end{displaymath}
  where $C$ is the norm of $c$ considered as an operator from $p_1H$ to $M$.

  2. We carry on the following computation, using now the fact that $c$ is a
  cocycle. For all $x\in\Ss$ and $y \in \Ss_1$ we have
  \begin{align*}
    m_c (\Lambda(x)\tens \Lambda(y)) &= xc(y) = c(xy) - c(x)\epsilon(y) \\
    &= c (xy - x\epsilon(y)) = c (\B-\O) (\Lambda(x)\tens \Lambda(y)),
  \end{align*}
  considering $c$ as a map from $\Ss \simeq \Hh$ to $M \subset H$.

  In particular this shows that $m_c$ vanishes on $\Ker(\B-\O) \cap \Kk$, which
  contains $K_g^\bot \cap \Kk$ by Lemma~\ref{lem_B-O}. Since $m_c$ is bounded
  and $K_g^\bot \cap \Kk$ is clearly dense in $K_g^\bot = \overline\Img
  (\Theta+\id)^*$, we can conclude that $m_c$ vanishes on $K_g^\bot$. As a result,
  the identity $m_c = c(\B-\O)$ holds on $\Kk \oplus K_g^\bot$, and in
  particular, on $\Kk'_g$.

  3. Now we apply this identity to the values of the path cocycle $c_g$, which
  lie in $\Kk'_g$ by definition. This yields, for any $x\in\Ss$:
  \begin{displaymath}
    m_c(c_g(x)) = c(\B-\O)(c_g(x)) = 2 c(x-\epsilon(x)1) = 2 c(x).
  \end{displaymath}
  Since $m_c$ is defined on the whole of $K$, we can also put $\xi = \frac 12
  m_c(\xi_g)$, and we have for any $x\in\Ss$:
  \begin{displaymath}
    x\xi - \xi\epsilon(x) = \ts\frac 12\, m_c(x\xi_g - \xi_g\epsilon(x)) = 
    \frac 12\, m_c(c_g(x)) = c(x),
  \end{displaymath}
  so that $\xi$ is a fixed vector in $H$ for $c : \Ss \to M$.

  4. It remains to check that the fixed vector $\xi$ lies in fact in $M$. But
  this is clear from the definition of $m_c$ and the hypothesis $\xi_g \in
  \Lambda(M)\tens p_1H$.
\end{dem}

The property of rapid decay gives a convenient way to check that a vector
$\xi\in H$ lies in $\Lambda_h(M)$, and it is known to hold for $A_o(Q)$ in the
unimodular case \cite[Thm.~4.9]{Vergnioux_RD}. Let us recall the definition of
Property RD.

We denote by $L$ be the classical length function on the discrete quantum group
under consideration, i.e. the central unbounded multiplier of $\hat S$ given by
$L = \sum l(\alpha) p_\alpha$, where $l(\alpha)$ is the distance from $1$ to
$\alpha$ in the classical Cayley graph. We consider $L$ as an unbounded operator
on $H$, with domain $\Hh$, and we denote by $H^s$ the completion of $\Ss$ with
respect to the norm $\|x\|_{2,s} = \|(1+L)^s\Lambda_h(x)\|_H$. Property RD
states that there exist constants $C>0$ and $s > 1$ such that $\|x\|_{S_\red}
:= \|\lambda(x)\|_{B(H)} \leq C \|x\|_{2,s}$ for all $x\in\Ss$. In other words
we have then continuous inclusions $H^s \subset S_\red \subset M \subset H$.

\begin{crl} \label{crl_vanishing_Ao} 
  Let $\Ss$, $M$ be the Hopf algebra and the von Neumann algebra associated with
  $A_o(I_n)$, $n\geq 3$. Then we have $H^1(\Ss, M) = 0$. 
\end{crl}

\begin{dem}
  According to the Theorem, it suffices to show that the fixed vector $\xi_g$ of
  the path cocycle lies in $\Lambda(M)\tens p_1H$. Thank to Property RD we will
  in fact prove that $\xi_g$ lies in $H^s\tens p_1H$, using
  Formula~\eqref{eq_fixed}.

  Since $(1+L)^s\tens\id = (k+1)^s\id$ on $p_k(K)$, the vectors
  $((1+L)^s\tens\id) (\tilde\xi_{i\wedge i+1})$ are pairwise orthogonal, and
  their respective norms are clearly dominated by $(i+2)^s$. As a result we have
  \begin{displaymath}
    \|\xi_g\|_{2,s}^2 \leq \frac 2{m_1} \sum_{i=0}^\infty \frac {(i+2)^{2s}}{m_im_{i+1}}.
  \end{displaymath}
  The left-hand side is finite because the sequence of quantum dimensions
  $(m_i)$ grows geometrically when $n\geq 3$, hence we are done.
\end{dem}

Let $(\Ss,\delta)$ be the dense Hopf algebra of a unimodular countable discrete
quantum group and denote by $M$ the von Neumann algebra associated with
$\Ss$. The definition of $L^2$-Betti numbers given by W.~L\"uck
\cite{Lueck_L2Betti} extends without difficulty to the quantum case
\cite{Kyed_Hom}, as well as the arguments of \cite{Thom_Cohom,Thom_Cocycles}, so
that we have
\begin{displaymath}
  \beta_k^{(2)}(\Ss,\delta) = \dim_M H_k(\Ss, M) = \dim_{M^\circ} H^k(\Ss, M).
\end{displaymath}
Here the action of $M^\circ$ comes from multiplication of elements of $M$ on the
right of the coefficient module $M$.

The result of Corollary~\ref{crl_vanishing_Ao} yields immediately the following
result, which strongly contrasts with the values $\beta_1^{(2)}(F_n) = n-1$:
\begin{crl}
  For any $n\geq 3$ we have $\beta_1^{(2)}(\Aa_o(I_n),\delta) = 0$.
\end{crl}
Notice that we also have $\beta_0^{(2)}(\Ss,\delta) = 0$ for $\Ss = \Aa_o(I_n)$,
because $M = A_o(I_n)''$ is diffuse. Moreover a free resolution of the co-unit
$\epsilon : \Ss \to \CC$ has been constructed in \cite{CHT_Homology_Ao}, showing
in particular that $\beta_k^{(2)}(\Ss,\delta) = 0$ for $k\geq 4$, and
$\beta_{3-k}^{(2)}(\Ss,\delta) = \beta_k^{(2)}(\Ss,\delta)$ for $k \in \{0, 1,
2, 3\}$. Using the Corollary above it results that $\beta_k^{(2)}(\Ss,\delta) =
0$ for all $k$. This also holds when $n = 1$, $2$ for a completely different
reason, namely amenability.

In the unitary case, Proposition~\ref{prp_unitary_nontriv} shows that $H^1(\Ss,
K) \neq 0$ for $\Ss = \Aa_u(I_n)$, $n\geq 1$. Since $K$ is a direct sum of
finitely many copies of $H$ as an $\Ss$-module, this implies $H^1(\Ss, H) \neq
0$. Recall that for non-amenable countable discrete groups, the vanishing of
$\beta_1^{(2)}(\Ss,\delta)$ is equivalent to the vanishing of $H^1(\Ss, H)$
(\cite{Bekka_Valette_Betti}, see also \cite[Cor.~2.4]{Thom_Cocycles}). This
carries over to discrete quantum groups \cite{Kyed_T}, hence we obtain
$\beta_1^{(2)}(\Aa_u(I_n),\delta) \neq 0$ for $n\geq 2$. It seems reasonable to
conjecture that $\beta_1^{(2)}(\Aa_u(I_n),\delta) = 1$.

\begin{rque}{Remark} {\sc A Hilbertian Variant.} \\
  We give now a direct proof of the vanishing of $H^1(\Aa_o(I_n), H)$ for $n
  \geq 3$, by adapting the methods of Theorem~\ref{thm_vanishing}.

  For a general cocycle $c : \Ss\to H$ one can define an unbounded map $m_c :
  \Kk \to H$ as in the Theorem, and it turns out that $m_c$ is closable as a
  densely defined operator on $K$. Consider indeed a sequence $\zeta_k =
  \sum_{ij} \Lambda(x_{ij,k}) \tens \Lambda(u_{ij}) \in \Kk$ that converges to
  $0$ and whose image under $m_c$ converges to some vector $\eta \in H$. We have
  then for any $z\in\Ss$:
  \begin{displaymath}
    \ts (\Lambda(z) | m_c(\zeta_k))
    = \sum (x_{ij,k}^* \Lambda(z) | c(u_{ij}))
    = \sum (\rho(z) \Lambda(x_{ij,k}^*) | c(u_{ij})),
  \end{displaymath}
  where $\rho(z)$, for $z\in \Ss$, is the bounded operator of right
  multiplication by $z$ in the GNS construction of $h$. We observe that the
  left-hand side converges to $(\Lambda(z)|\eta)$, whereas the right-hand side
  converges to $0$ since the sum is finite and the convergence $\zeta_k \to 0$
  means that $\Lambda(x_{ij,k}) \to 0$ for each $i$, $j$. This shows that $\eta
  = 0$.

  Now the vanishing of $m_c$ on $\Ker(\B-\O) \cap \Kk$ implies the vanishing of
  its closure $\overline m_c$ on $\Ker(\B-\O) \cap \Dom \overline m_c$, so that
  the identity $\overline m_c = c \rond (\B-\O)$ holds on $\Kk'_g \cap \Dom
  \overline m_c$. But it is easy to check that $\Lambda_h(M)\tens p_1H \subset
  \Dom \overline m_c$, and since the fixed vector $\xi_g$ lies in
  $\Lambda_h(M)\tens p_1H$ we obtain the inclusion $c_g(\Ss) \subset \Dom
  \overline m_c$. We can then apply $m_c$ to $c_g(x)$ and $\xi_g$ as in the
  proof of the Theorem and obtain the triviality of $c$.
\end{rque}

\begin{rque}{Remark} {\sc The non-unimodular case.} \\
  $L^2$-Betti numbers have not been defined for discrete quantum groups which
  are not unimodular. However we can still prove the vanishing of an
  $L^2$-cohomology group in this case, i.e. for $A_o(Q)$ with $n\geq 3$, $Q\bar
  Q \in \CC I_n$, $Q$ not unitary.

  More precisely it is natural for the proof of Theorem~\ref{thm_vanishing} to
  consider the submodule $N\subset H$ of vectors $\xi$ such that $\rho(\xi) :
  \Hh \to H$, $\Lambda(x)\mapsto x\xi$ is bounded. The left $\Ss$-module $N$ can
  be identified with $M$ endowed with the twisted left action $a \cdot x =
  \sigma^h_{i/2}(a)x$, via $(x \mapsto Jx^*\xi_0)$. In other words we are
  considering the group $H^1(\Ss,M)$ with the twisted action on the left of $M$,
  and the trivial action on the right. Let us explain how the techniques of
  Theorem~\ref{thm_vanishing} can be used to prove the vanishing of this group.

  For any cocycle $c : \Ss \to N$ one can define a bounded map $m_c : K \to H$
  by the formula $m_c(\xi\tens\Lambda(y)) = \rho(c(y))\xi$ and follow the proof
  of the Theorem. The only point that requires more care is the fact that
  $\xi_g$ lies in $N\tens p_1H$ --- since $\rho(\Ss)N \subset N$, this will
  imply that the fixed vector $\xi$ belongs to $N$. We cannot apply Property RD
  anymore, since it does not hold for non-unimodular discrete quantum
  groups. However, we have the following strong weakening of Property RD,
  already used in \cite[Rem.~7.6]{VaesVer_Boundary}, which will suffice for our
  purposes.

  Denoting by $H^s_r$ the completion of $\Hh$ with respect to the norm
  $\|\zeta\|_{2,s,r} = \|r^L (1+L)^s\zeta\|_H$, the methods of
  \cite{Vergnioux_RD} yield, in the non-unimodular case, a continuous inclusion
  $H^s_r \subset \Lambda(M)$, where $s=3$, $r = \|F\|$ and $F$ is the usual
  positive ``modular'' element of $B(H_1)$ satisfying $\Tr (F) = \Tr
  (F^{-1})$. Note that we have also $H^s_r \subset N$ since $J H^s_r = H^s_r$
  --- indeed $r^L (1+L)^s$ is scalar on each $p_k H$ and $J : p_k H \to p_kH$ is
  an (anti-)isometry.

  Therefore it suffices to check that $\xi_g$ lies in $H^s_r \tens p_1H$. Since
  $L = k\,\id$ on $p_k H$ we get immediately from~\eqref{eq_fixed}:
  \begin{displaymath}
    \|\xi_g\|_{2,s,r}^2 \leq \frac 2{m_1}
    \sum_{i=0}^\infty \frac {r^{2i+2}(i+2)^{2s}}{m_i m_{i+1}}
    \leq C \sum_{i=0}^\infty \big(\ts\frac ra\big)^{2i} (i+2)^{2s},
  \end{displaymath}
  where $a> 1$ is the number such that $\Tr(F) = a + a^{-1}$ --- we have then
  $m_k = (a^{k+1}-a^{-k-1})/(a-a^{-1})$ for all $k$, hence the second upper
  bound above. Now we have $a + a^{-1} = \Tr(F) > r+r^{-1}$ since $r$, hence
  also $r^{-1}$, is an eigenvalue of $F$, and the size of $F$ is at least
  $3$. This implies $a > r$, hence $\|\xi_g\|_{2,s,r} < \infty$.
\end{rque}


\end{document}